\documentclass[a4paper,12pt]{amsart}
\usepackage[T1]{fontenc} % codifica dei font
\usepackage[utf8]{inputenc} % lettere accentate da tastiera
\usepackage[english]{babel} % lingua del documento
\usepackage{amsmath}
\usepackage{amssymb}
\usepackage[english]{varioref}
\usepackage{amsthm}
\usepackage{mathrsfs}
\usepackage[all,cmtip]{xy}
\usepackage{enumitem}
\usepackage{graphicx}
\usepackage{caption}
\usepackage{hyperref}
\usepackage{longtable}

\hypersetup{
	colorlinks=true,
	linkcolor=blue,}

\calclayout

\let\OLDthebibliography\thebibliography
\renewcommand\thebibliography[1]{
	\OLDthebibliography{#1}
	\setlength{\parskip}{0pt}
	\setlength{\itemsep}{5pt}
}

\newcommand{\ze}[1]{\zeta_{#1}}
\def\t{\text}

\DeclareMathOperator{\Z}{\mathbb{Z}}
\DeclareMathOperator{\Q}{\mathbb{Q}}
\DeclareMathOperator{\HOM}{\mathrm{Hom}}
\DeclareMathOperator{\A}{\mathcal{A}}
\DeclareMathOperator{\cyc}{\mathrm{Cyc}}
\DeclareMathOperator{\norm}{\mathrm{Norm}}
\DeclareMathOperator{\gal}{\mathrm{Gal}}
\DeclareMathOperator{\cl}{\mathrm{Cl}}
\DeclareMathOperator{\ext}{\mathrm{Ext}^1}

\newtheorem{theorem}{Theorem}[section]
\newtheorem*{theorem*}{Theorem}
\newtheorem{proposition}[theorem]{Proposition}
\newtheorem{lemma}[theorem]{Lemma}

\newtheorem{corollary}[theorem]{Corollary}

\theoremstyle{definition}
\newtheorem{definition}[theorem]{Definition}
\newtheorem{remark}[theorem]{Remark}

\begin{document}

\title[Greenberg's conjecture for real quadratic fields]{Greenberg's conjecture for real quadratic fields and the cyclotomic $\Z_2$-extensions}
	
\author{Lorenzo Pagani}
\address{Dipartimento di matematica, università di Roma ``La Sapienza'', piazzale Aldo Moro 5, 00185 Roma Italy}

\email{lorenzo.pagani3110@gmail.com}

%%\date{\today}

\begin{abstract}
	Let $\A_n$ be the $2$-part of the ideal class group of the $n$-th layer of the cyclotomic $\Z_2$-extension of a real quadratic number field $F$. The cardinality of $\A_n$ is related to the index of cyclotomic units in the full group of units. We present a method to study the latter index. As an application we show that the sequence of the $\A_n$'s stabilizes for the real fields $F=\Q(\sqrt{f})$ for any integer $0<f<10000$. Equivalently Greenberg's conjecture holds for those fields.    
\end{abstract}

\maketitle

\section*{Introduction}

Let $F$ be a number field. A $\Z_p$-extension is a field extension $\mathcal{F}|F$ whose Galois group is topologically isomorphic to the additive group of $p$-adic numbers. Moreover, $\mathcal{F}$ is said to be cyclotomic if it is the compositum between $F$ and the unique $\Z_p$-extension of the rationals. 
We have a tower of fields
\[ F=F_0\subset F_1\subset F_2\subset \dots \subset \mathcal{F} \]
such that $\gal(F_n|F)$ is a cyclic group of order $p^n$. The $p$-parts $\A_n$ of the class groups of the fields $F_n$ give a projective system
\[ \A_0\overset{N}{\longleftarrow}\A_1\overset{N}{\longleftarrow}\A_2\overset{N}{\longleftarrow}\dots \]
where $N$ is the norm map. In \cite{iwasawa}, Iwasawa showed the existence of integers $\mu$, $\lambda$ and $\nu$ such that
\[ \#\A_n=p^{\mu p^n+\lambda n+\nu} \]
for large enough $n$. 

For cyclotomic $\Z_p$-extensions of abelian fields $F$, it is known that $\mu=0$ as shown by Ferrero and Washington in \cite{ferrerowash}. Moreover, Greenberg in his thesis \cite{green} conjectured that $\lambda=0$ when $F$ is a real field. In other words, Greenberg's conjecture is equivalent to say that the cardinalities of $\A_n$ stabilize in the cyclotomic tower. 

In this paper we restrict our attention to the case in which $F$ is a real quadratic field and $\mathcal{F}$ is the cyclotomic $\Z_2$-extension of $F$. This problem has been studied by various authors. In \cite{ozakitaya}, \cite{mizusawa}, \cite{nishino}, \cite{fukudakomatsuprev} and \cite{fukudakomatsu} Greenberg's conjecture has been proven when $F=\Q(\sqrt{f})$ where $f$ is divisible by at most $3$ primes satisfying some arithmetic properties.

We follow a different approach. We avoid studying directly the module $\A_n$. Instead, we study cyclotomic units in order to define modules $A_n$ whose cardinality is closely related to the one of $\A_n$, then exploiting properties of Gorenstein rings we give a procedure to compute the modules $A_n$. This method is an adaptation of the algorithm described in \cite{corn}, \cite{kraftschoof} or \cite{schoof2} to our situation. As a consequence of our computation we have the following.
\begin{theorem*}
	Greenberg's conjecture is true for the cyclotomic $\Z_2$-extension of real quadratic fields $F=\Q(\sqrt{f})$ where $f$ is an integer and $0<f<10000$.
\end{theorem*}
Similar computational results for cyclotomic $\Z_3$-extensions of quadratic real fields were already known. They can be found in \cite{fukuda}, \cite{taya} and \cite{kraftschoof}. We follow the approach of the last article. 

In adapting it to our situation we encounter two main obstacles: the fact that for the field $F_n$ the index of cyclotomic units inside the full group of units   
is equal to the class number $\#\cl_{F_n}$ times a $2$-power factor that grows as $2^{[F_n:\Q]}$. The second one is the fact that when $F$ is a quadratic field, a $\Z_{2}[\gal(F|\Q)]$-module in general does not split as a part which is invariant with the action of $\gal(F|\Q)$ and a part which is anti-invariant. While it is true for $\Z_p[\gal(F|\Q)]$-modules for odd primes $p$. Furthermore, we expect that the algorithm we describe can be extended to  the case of cyclotomic $\Z_p$-extensions of Galois extensions of $\Q$ of degree $p$, where $p$ is an odd prime integer. We think that this case could be easier since $1$ is the only real root of unity of order a power of $p$ when $p$ is odd. 

The results presented in \cite{kraftschoof} has already been extended in different directions. Paoluzi in \cite{paoluzi} studied the case in which the odd prime $p=3$ splits in the quadratic field. This case has been handled in a way that is different from the approach presented in this paper. Moreover, Mercuri extended the computations also to the cases $p=5$ and $p=7$ in \cite{mercuri}. 

In the first section we study the modules $A_n$ defined as the $2$-parts of the quotient of $O_{F_n}^{\times}$ by a suitable group of units $C_{F_n}$. The group $C_{F_n}$ is closely related to the group of Sinnott cyclotomic units, this relation is made explicit in the appendices. In particular, the cardinality of $A_n$ is equal to the cardinality of the $2$-part of $\#\cl_{F_n}$ or half of it according to the behaviour of the prime $2$ in $F$. Thus, proving that the modules $A_n$ stabilize is equivalent to proving Greenberg's conjecture for a fixed quadratic field $F$. In section $2$ and $3$ we exploit properties of finite Gorenstein rings in order to describe the $\gal(F_n|F)$-module structure of the dual module of $A_n$. Finally, in section $4$ we present some criteria that permit us to deduce that the modules $A_n$ stabilize and therefore Greenberg's conjecture holds. In section $5$ we present some examples, making some relevant remarks on the algorithm, and we collect the obtained results in section $6$.

\section{The quotient of units modulo cyclotomic units}	

Let $f\ge3$ be a squarefree integer, let $F=\Q(\sqrt{f})$. Denote by $F_n$ and $\Q_n$ the $n$-th levels of the cyclotomic $\Z_2$-extensions of $F$ and $\Q$ respectively. In other words, we have $\Q_n=\Q(\ze{2^{n+2}})^+$ and $F_n=\Q_n(\sqrt{f})$ where $\ze{2^{n+2}}$ is a primitive root of unity of order $2^{n+2}$. We exclude the case $f=2$ since $\Q_1=\Q(\sqrt{2})$ and therefore the cyclotomic $\Z_2$-extensions of $\Q$ and $\Q(\sqrt{2})$ coincide. Recall that the class number of $\Q_n$ is odd, this well known fact can be deduced from Chevalley's formula following to the idea explained in Remark \ref{2-principality} presented below. Similarly, we can assume $f$ to be odd. Indeed, since $\Q_1=\Q(\sqrt{2})$ we have that the cyclotomic $\Z_2$-extension of $\Q(\sqrt{f})$ and $\Q(\sqrt{2f})$ coincide for all $n>0$. 

We fix the following notation: let $\Gamma_n=\gal(F_n|\Q)$. We have that $\Gamma_n\cong G_n\times\Delta$, where $G_n=\gal(F_n|F)\cong\gal(\Q_n|\Q)$ is cyclic of order $2^n$ while $\Delta=\gal(F_n|\Q_n)\cong\gal(F|\Q)$ is cyclic of order $2$.

Denote by $\A_n$ the $2$-part of the class group of $F_n$. We have that the norm maps $\A_{n+1}\longrightarrow \A_n$ are surjective for any $n$ since $F_{n+1}$ and the Hilbert's class field of $F_n$ are linearly disjointed. The Chevalley's formula shown in \cite{chevalley} gives us a relation between the cardinalities of the groups $\A_n$ inside the $\Z_p$-tower. We have that for any $0\le m\le n$
\begin{equation}\label{Chev-formula} \tag*{(C)}
	 \#\A_n^{\gal(F_n|F_m)} = \begin{cases} \#\A_m, & \t{if $2$ does not split in $F$,}
\\ 2^h\#\A_m, & \t{if $2$ splits in $F$}
\end{cases} 
\end{equation}
where $0\le h\le n-m$. The factor $2^h$ is a normic factor of units obtained from the fundamental product formula of class field theory. The definition of this factor can be found in \cite{chevalley} or \cite[Remarks 6.2.3]{gras}. Furthermore, in \cite{fukuda-stab} it is shown that if there exists a level $n_0$ such $\#\A_{n_0}=\#\A_{n_0+1}$, then the $\lambda$-invariant is zero. However, we do not have a proof of the existence of such $n_0$.

\begin{remark}\label{2-principality}
	Assume that $2$ does not split in $F$. Then, by Chevalley's  formula \ref{Chev-formula} we get that $\#\A_n^{G_n}=\#\A_0$. Since both $G_n$ and $\A_n$ are $2$-groups we deduce that $\A_n$ is trivial if and only if $\A_0$ is. 
\end{remark}

So, we can make the general assumption that $f\ge3$ is an odd squarefree integer. Moreover if $2$ does not split in $F=\Q(\sqrt{f})$, we assume that $2|\#\cl_{F}$. Note that $2$ does not split in $F$ whenever $f\not\equiv1\bmod8$. Finally, we recall that $f$ is not necessarily the conductor of the field $F$.

We define a module $C_{F_n}$ contained in $O_{F_n}^{\times}$ whose index is closely related to the class number of the field of $F_n$. Then, we denote by $A_n$ the $2$-part of the quotient $O_{F_n}^{\times}/C_{F_n}$ and we study $A_n$ instead of the $2$-part of $\cl_{F_n}$.  We will show that there are some analogies between $\A_n$ and $A_n$, for instance compare Lemma \ref{relation-An} with Chevalley's formula and Proposition \ref{An-stabilizatio} with \cite[Theorem 1]{fukuda-stab}.

\begin{definition}\label{moreunits}
	Fix an algebraic closure $\overline{\Q}$ of the rationals. For any integer $n\ge0$ define the following elements:
	\begin{enumerate}[label=(\roman*)]
		\item Let $\sigma$ be the generator of $\Delta$. For any positive squarefree integer $f$ such that $f\equiv1\bmod4$ let $\delta_f\in\overline{\Q}$ be such that
		\[ \delta_f^2=\begin{cases}
		\left(\underset{\Q(\ze{f})|F}{\norm}(1-\ze{f})\right)^{1-\sigma}, & \t{if $f$ prime,}\\
		\underset{\Q(\ze{f})|F}{\norm}(1-\ze{f}), & \t{otherwise.}
		\end{cases} \]
		Note that $\delta_f$ is defined up to a sign and $\delta_f^2$ lies in $O_{F_0}^{\times}$.
		\item Let $\gamma$ be a generator of $G_n\cong\gal(\Q_n|\Q)$. Fix an element $c\in(\Z/2^{n+2}\Z)^\times$ such that $\gamma(\ze{2^{n+2}}+\ze{2^{n+2}}^{-1})=\ze{2^{n+2}}^c+\ze{2^{n+2}}^{-c}$. Then, define  
		\[ \beta_n=\ze{2^{n+2}}^{\frac{1-c}{2}}\frac{1-\ze{2^{n+2}}^c}{1-\ze{2^{n+2}}}\in O_{\Q_n}^\times. \]
		Note that $c$ is defined up to a sign. Choosing $-c$ instead of $c$ changes the sign of the unit $\beta_n$. 
		\item We define the element of $F_{n+1}^\times$
		\[ \eta_n=\underset{\Q(\ze{2^{n+3}f})^+|F_{n+1}}{\norm}\left( \ze{4}(\ze{2^{n+3}f}-\ze{2^{n+3}f}^{-1})\right).  \]
	\end{enumerate}
\end{definition}

These elements are related to the cyclotomic units introduced by Sinnott in \cite[Section $4$]{sinnott}. We make this relation explicit in Lemma \ref{cyclotomic-units}.

\begin{definition}\label{new-group-of-cyclounits}
	Let $f\ge3$ be an odd squarefree integer. Assume that the class number of $\Q(\sqrt{f})$ is even if $f\not\equiv1\bmod8$. Denote by $\delta_f'$ the unit $\delta_f$ if $f\equiv1\bmod 8$ and $f$ is prime or the unit $\delta_f^2$ if $f\equiv1\bmod8$ and $f$ is not prime. For any $n\ge0$ we define the following $G_n$-modules $C_{F_n}\subset O_{F_n}^{\times}$:
	\[ C_{F_n}=\begin{cases}\left\langle -1,\eta_n,\beta_n,\delta_f'\right\rangle_{\Z[G_n]}, & \t{if $f\equiv1\bmod8$ }, \\
	\left\langle -1,\eta_n,\beta_n\right\rangle_{\Z[G_n]}, & \t{if $f\not\equiv1\bmod8$ }.
	 \end{cases} \] 
\end{definition} 

The fact that $C_{F_n}$ is contained in $O_{F_n}^{\times}$ is proven in Appendix \ref{appB} as well as the following theorem:

\begin{theorem}\label{CFn-structure}
	Let $f\ge3$ be an odd squarefree integer. Assume that the class number of $\Q(\sqrt{f})$ is even if $f\not\equiv1\bmod8$. The modules $C_{F_n}$ satisfy the following properties:
	\begin{enumerate}[label=(\alph*)]
		\item for any $n\ge0$, the index $[O_{F_n}^{\times}:C_{F_n}]$ is equal to $\#\cl_{F_n}/2$ if $f\equiv3\bmod4$ or is equal to $\#\cl_{F_n}$ otherwise.
		\item\label{gal-prop-CFn} Let $0\le m\le n$. Then, we have that $C_{F_m}=C_{F_n}^{\gal(F_n|F_m)}$.
	\end{enumerate}
	Let $\gamma$ be a generator of $G_n=\gal(F_n|F)$, define $N_n$ as the element  $\sum_{i=0}^{2^n-1}\gamma^i\in\Z[G_n]$. Assume that $f\equiv1\bmod8$. Then, for any $n\ge0$ we have an isomorphism
	\[ \pi:\frac{\Z[G_n]}{(N_n)}\oplus\frac{\Z[G_n]}{(N_n)}\oplus\Z\longrightarrow\frac{C_{F_n}}{\pm1} \]
	defined by $(x,y,z)\mapsto\eta_n^x\beta_n^y{\delta_f'}^z$. Here $\delta_f'=\delta_f$ if $f$ is prime or $\delta_f'=\delta_f^2$ otherwise. 
	
	Assume now that $f\not\equiv1\bmod8$. Then, for any $n\ge0$ we have an isomorphism
	\[ \pi:\Z[G_n]\oplus\frac{\Z[G_n]}{(N_n)}\longrightarrow\frac{C_{F_n}}{\pm1} \]
	defined by $(x,y)\mapsto\eta_n^x\beta_n^y$.
\end{theorem} 

\begin{definition}\label{An}
	Let $f\ge3$ be an odd squarefree integer. Assume that the class number of $\Q(\sqrt{f})$ is even if $f\not\equiv1\bmod8$. For any $n\ge0$ we define the $\Z_2[G_n]$-module
	\[ A_n=\frac{O_{F_n}^{\times}}{C_{F_n}}\otimes_{\Z}\Z_2. \]
\end{definition}
\begin{remark}\label{An-injectivity}
	Note that the cardinality of $A_n$ is closely related to the $2$-part of the class number of $F_n$. For any $0\le m\le n$ the inclusions $O_{F_m}^{\times}\subset O_{F_n}^{\times}$ induce natural morphisms $A_m\longrightarrow A_n$. By Theorem \ref{CFn-structure} part \ref{gal-prop-CFn} these morphisms are injective. In particular, the sequence of the cardinalities of $A_n$ is non-decreasing. Proving that this sequence stabilizes is equivalent to proving Greenberg's conjecture for the field $F$.  
\end{remark}

\begin{lemma}\label{relation-An}
	Let $f\ge3$ be an odd squarefree integer. Assume that the class number of $\Q(\sqrt{f})$ is even if $f\not\equiv1\bmod8$. For any $0\le m\le n$ let $\gamma$ be a generator of $G_n$ and let $N_m^{(n)}$ be the element of $\Z[G_n]$ defined as $\sum_{i=0}^{2^m-1}\gamma^i$. Then, the natural map $A_m\longrightarrow A_n$ induces an isomorphism
	\[ \frac{A_m}{A_0}\cong\frac{A_n}{A_0}\left[N_m^{(n)}\right]. \]
	Here $(A_n/A_0)[N_m^{(n)}]$ denotes the kernel of the map $A_n/A_0\longrightarrow A_n/A_0$ induced by the multiplication by $N_m^{(n)}$.
\end{lemma}

\begin{proof}
	We consider the following commutative diagram with exact rows:
	\[\xymatrix{0\ar[r]&C_{F_0}\ar[r]\ar[d]&O_{F_0}^\times\ar[r]\ar[d]&\frac{O_{F_0}^\times}{C_{F_0}}\ar[r]\ar[d]&0
		\\0\ar[r]&C_{F_n}\ar[r]&O_{F_n}^\times\ar[r]&\frac{O_{F_n}^\times}{C_{F_n}}\ar[r]&0 }\]
	By Theorem \ref{CFn-structure} part \ref{gal-prop-CFn} we have that $C_{F_0}=C_{F_n}^{G_n}$. Hence, the rightmost vertical map is injective. Therefore, the snake lemma gives us the following exact sequence:
	\[ 0\longrightarrow \frac{C_{F_n}}{C_{F_0}}\longrightarrow \frac{O_{F_n}^\times}{O_{F_0}^\times}\longrightarrow \frac{O_{F_n}^\times}{C_{F_n}O_{F_0}^\times}\longrightarrow0. \]
	The cardinality of the rightmost group of the exact sequence is $\#\cl_{F_n}/\#\cl_{F_0}$.
	
	According to Definition \ref{new-group-of-cyclounits} we have that $C_{F_0}=\left\langle -1,\delta_f'\right\rangle_{\Z}$ if $f\equiv1\bmod8$ or we have that $C_{F_0}=\left\langle -1,\eta_0\right\rangle_{\Z} $. Indeed, the unit $\beta_0$ is equal to $-1$ and by Remark \ref{norm-relation} we get that up to a sign
	\[ \eta_0=\begin{cases}1, & \t{if $f\equiv1\bmod8$},\\
	\underset{F_n|F_0}{\norm}(\eta_n), & \t{otherwise}.	
	\end{cases} \]
	By Theorem \ref{CFn-structure} we deduce that the group $C_{F_n}/C_{F_0}$ is isomorphic as a $\Z[G_n]$-module to $\Z[G_n]/(N_n)\oplus\Z[G_n]/(N_n)$.

	We apply the functor $\HOM_{\Z[G_n]}(\Z[G_n]/(N_m^{(n)}),-)$ to get
	\[ 0\longrightarrow \frac{C_{F_n}}{C_{F_0}}\left[N_m^{(n)}\right]\longrightarrow \frac{O_{F_n}^\times}{O_{F_0}^\times}\left[N_m^{(n)}\right]\longrightarrow \frac{O_{F_n}^\times}{C_{F_n}O_{F_0}^\times}\left[N_m^{(n)}\right]\longrightarrow0. \]
	The above short sequence is exact since $\ext_{\Z[G_n]}(\Z[G_n]/(N_m^{(n)}),\Z[G_n]/(N_n))=0$. This can be deduced by the fact that $N_m^{(n)}$ divides $N_n$ in the module $\Z[G_n]$ for any $0\le m\le n$. 
	
	Furthermore, we have the following commutative diagram whose rows are exact:
	\[\xymatrix{0\ar[r]&\frac{C_{F_m}}{C_{F_0}}\ar[r]\ar[d]&\frac{O_{F_m}^\times}{O_{F_0}^\times}\ar[r]\ar[d]&\frac{O_{F_m}^\times}{C_{F_m}O_{F_0}^\times}\ar[r]\ar[d]&0
		\\0\ar[r]&\frac{C_{F_n}}{C_{F_0}}\left[N_m^{(n)}\right]\ar[r]&\frac{O_{F_n}^\times}{O_{F_0}^\times}\left[N_m^{(n)}\right]\ar[r]&\frac{O_{F_n}^\times}{C_{F_n}O_{F_0}^\times}\left[N_m^{(n)}\right]\ar[r]&0 }.\]
	The leftmost vertical map is an isomorphism. Indeed, the map
	\[ \frac{C_{F_m}}{C_{F_0}}\longrightarrow\frac{C_{F_n}}{C_{F_0}}\left[N_m^{(n)}\right] \] 
	is clearly injective. Take $x\in C_{F_n}$ such that $x^{N_m^{(n)}}\in C_{F_0}$. Then $1=\left(x^{N_m^{(n)}}\right)^{\gamma-1}=x^{\gamma^{2^m}-1}$. Therefore $x\in C_{F_m}$ and the surjectivity follows. Similarly, we prove that also the central vertical map is an isomorphism.  
	Hence, we get that
	\[ \frac{O_{F_m}^\times}{C_{F_m}O_{F_0}^\times}\cong \frac{O_{F_n}^\times}{C_{F_n}O_{F_0}^\times}\left[N_m^{(n)}\right]. \]
	
	Since we are interested in the $2$-part of the class group, we tensor by $\Z_2$ and we get
	\[ \frac{A_m}{A_0}\cong \frac{A_n}{A_0}\left[N_m^{(n)}\right]. \]
\end{proof}

\begin{proposition}\label{An-stabilizatio}
	Let $f\ge3$ be an odd squarefree integer. Assume that the class number of $\Q(\sqrt{f})$ is even if $f\not\equiv1\bmod8$. Assume there exists an integer $n_0\ge0$ such that
	$A_{n_0}$ and $A_{n_0+1}$ have the same cardinality. Then, for any $n\ge n_0$ the natural maps
	$A_{n_0}\longrightarrow A_n$ are isomorphisms.
\end{proposition}

\begin{proof}
	Fix $n>n_0$. By Lemma \ref{relation-An} we have
	\[ \frac{A_{n_0}}{A_0}\cong\frac{A_n}{A_0}\left[N_{n_0}^{(n)}\right]\subset\frac{A_n}{A_0}\left[N_{n_0+1}^{(n)}\right]\cong\frac{A_{n_0+1}}{A_0}. \]
	The hypothesis $\#A_{n_0}=\#A_{n_0+1}$ implies that $\#(A_{n_0}/A_0)=\#(A_{n_0+1}/A_0)$ too. Indeed, the natural maps $A_0\longrightarrow A_n$ are injective by the fact that $C_{F_0}=C_{F_n}^{G_n}$. Since all the modules are finite we get
	\[ \frac{A_n}{A_0}\left[N_{n_0}^{(n)}\right]=\frac{A_n}{A_0}\left[N_{n_0+1}^{(n)}\right]\quad \t{and} \quad N_{n_0}^{(n)}\frac{A_n}{A_0}=N_{n_0+1}^{(n)}\frac{A_n}{A_0}.\]
	
	The algebra $\Z_2[G_n]$ is local with maximal ideal $\mathfrak{M}=(2,\gamma-1)$. Note that
	\[ N_{n_0}^{(n)}\frac{A_n}{A_0}=N_{n_0+1}^{(n)}\frac{A_n}{A_0}=(1+\gamma^{2^{n_0}})N_{n_0}^{(n)}\frac{A_n}{A_0}\subset \mathfrak{M}N_{n_0}^{(n)}\frac{A_n}{A_0}. \]
	Therefore Nakayama's lemma implies that $N_{n_0}^{(n)}(A_n/A_0)=0$. Hence
	\[ \frac{A_{n_0}}{A_0}\cong\frac{A_n}{A_0}\left[N_{n_0}^{(n)}\right]=\frac{A_{n}}{A_0} . \]
	This implies that $\#A_{n_0}=\#A_n$.  Since the natural morphism $A_{n_0}\longrightarrow A_n$ is injective as stated in Remark \ref{An-injectivity} we deduce that it is an isomorphism. 
\end{proof}

This result about the stabilization of $A_n$ can also be seen as a consequence of \cite[Theorem 1]{fukuda-stab}. Interestingly, we are able to proving it using cyclotomic units.

Finally, in order to understand the structure of $ A_n$ we study its $2^k$-torsion part for a fixed positive integer $k$. If $k$ is big enough we have that $ A_n= A_n[2^k]$.

Consider the commutative diagram with exact rows:
\[ \xymatrix{0\ar[r]&\frac{C_{F_n}}{\pm1}\ar[r]\ar[d]^-{2^k}&\frac{O_{F_n}^\times}{\pm1}\ar[r]\ar[d]^-{2^k}&\frac{O_{F_n}^\times}{C_{F_n}}\ar[r]\ar[d]^-{2^k}&0
	\\0\ar[r]&\frac{C_{F_n}}{\pm1}\ar[r]&\frac{O_{F_n}^\times}{\pm1}\ar[r]&\frac{O_{F_n}^\times}{C_{F_n}}\ar[r]&0.}  \]  
The snake lemma gives us the exact sequence of $\Z[G_n]$-modules
\begin{equation}\label{torsion}\tag*{($1.1$)}
0\longrightarrow  A_n[2^k]\longrightarrow \frac{C_{F_n}}{\pm C_{F_n}^{2^k}}\longrightarrow \frac{C_{F_n}}{\pm O_{F_n}^{\times2^k}\cap C_{F_n}}\longrightarrow0.
\end{equation}

We will see in the next section that it is convenient to deal with the dual sequence.

\section{Dual modules over finite Gorenstein rings} 

In this section we recall the definition and some properties of finite Gorenstein rings. We use them to describe the dual of $C_{F_n}/\pm O_{F_n}^{\times 2^k}\cap C_{F_n}$ in terms of Galois groups. 

\begin{definition}
	Let $R$ be a commutative finite ring. For any $R$-module $A$, we have the group $A^{\perp}=\HOM_R(A,R)$ and the dual group $A^{\mathrm{dual}}=\HOM_{\Z}(A,\Q/\Z)$. They are $R$-modules via $(r\phi)(a)=\phi(ra)$ for any  $r\in R$ and $a\in A$.
	The ring $R$ is said to be Gorenstein if $R^{\mathrm{dual}}$ is a free $R$-module of rank one.
\end{definition}

We state below the properties of Gorenstein rings we need. 
\begin{proposition}\label{Gore-propeties}
	Let $R$ be a finite Gorenstein ring. Then
	\begin{enumerate}[label=(\alph*)]
		\item\label{a} For every $R$-module $A$, the map
		\[A^\perp=\HOM_R(A,R)\rightarrow A^{\mathrm{dual}}=\HOM_{\Z}(A,\Q/\Z)\]
		defined by $\psi\mapsto \phi\circ \psi$ is an isomorphism of $R$-modules. Here $\phi:R\rightarrow \Q/\Z$ denotes a generator of $R^{\mathrm{dual}}$.
		In particular if $A$ is finite, then $\#A=\#A^{\perp}$.
		\item The functor $A\rightarrow A^\perp$ from the category of finite $R$-module to itself is exact. 
	\end{enumerate}
\end{proposition} 
These properties are shown in \cite[Proposition $1.1$]{schoof2}.
\begin{remark}\label{Gore-isomorphism}
	Let $R_k=\Z/2^k\Z[G_n]$. It is a finite Gorenstein ring. We make explicit the isomorphism described in the proposition above. Let $\phi:R_k\longrightarrow\Z/2^k\Z$ be the homomorphism that maps an element $\sum_{g\in G_n}x_gg$ to $x_1$. Then $\phi$ generates the $R_k$-module $\HOM_{R_k}(R_k,\Z/2^k\Z)$ which is isomorphic to $R_k^{\mathrm{dual}}$. If $\psi$ is an element of $\HOM_{R_k}(A,R_k)$ that maps $a$ to $\sum_{g\in G_n}x_gg$, then the composition $\phi\circ\psi:A\longrightarrow\Z/2^k\Z$ maps $a$ to $x_1$.
\end{remark}

\begin{lemma}\label{dual-kill}
	Let $R$ be a finite Gorenstein ring and let $A$ be an $R$-module. Assume there exists $r\in R$ such that $r(A^{\perp})=0$. Then, the module $rA$ is trivial. 
\end{lemma}
\begin{proof}
	Consider the exact sequence of $R$-modules
	\[ 0\longrightarrow A[r]\longrightarrow A\overset{r}{\longrightarrow}A\longrightarrow A/rA\longrightarrow0. \]
	The dual map $A^\perp\longrightarrow A^{\perp}$ induced by the multiplication by $r$ on $A$ is the multiplication by $r$ on $A^\perp$. Therefore, applying the functor $(-)^\perp$ we get the exact sequence 
	\[0\longrightarrow (A/rA)^\perp\longrightarrow A^\perp\overset{r}{\longrightarrow}A^\perp\longrightarrow A[r]^\perp\longrightarrow0.\]
	Since $r$ annihilates $A^{\perp}$, we get $(A/rA)^{\perp}\cong A^\perp$. By the faithfulness of $(-)^\perp$, we deduce that the natural quotient $A\longrightarrow A/rA$ is an isomorphism. This implies that $rA=0$. 
\end{proof}

The aim of this section is describing the elements of $\left(C_{F_n}/\pm O_{F_n}^{\times 2^k}\cap C_{F_n} \right)^{\perp}$. 

\begin{definition}\label{function_fr}
	Fix $f\ge3$ an odd squarefree integer. Let $r$ be a prime that completely splits in $F_n(\ze{4})=\Q(\ze{2^{n+2}},\sqrt{f})$. In other words, $r\equiv1\bmod 2^{n+2}$ and $f$ is a square modulo~$r$. Let $k\le n+1$ be a positive integer. Choose a prime $\mathfrak{R}$ of $F_n(\ze{4})$ dividing $r$ and a primitive $2^k$-root of unity $\ze{2^k}$. For any $u\in O_{F_n}^{\times}$ we define $\log_r(u)$ as the element of $\Z/2^k\Z$ such that
	\[ u^{\frac{r-1}{2^k}}\equiv\ze{2^k}^{\log_r(u)} \bmod \mathfrak{R}. \]
	Since $k\le n+1$, the function $\log_r$ is trivial when evaluated on elements $u\in\pm O_{F_n}^{\times 2^k}$.
	We define $f_r\in\left(C_{F_n}/\pm O_{F_n}^{\times 2^k}\cap C_{F_n} \right)^{\perp}$ as
	\[ f_r(u)=\sum_{g\in G_n}\log_r(u^g)g^{-1}. \]
\end{definition}

Notice that the definition of $f_r$ depends both on a choice of a primitive $2^k$-th root of unity and of a prime $\mathfrak{R}$. If we choose a different $\ze{2^k}$ we are multiplying $f_r$ by an invertible element of $\Z/2^k\Z$. While if we change $\mathfrak{R}$ we are multiplying $f_r$ by an element of $G_n$. Therefore $f_r$ is well defined up to an invertible element of $R_k$. If we want to stress the choice of the ideal $\mathfrak{R}$ we write $f_{\mathfrak{R}}$ instead of $f_r$. 

\begin{proposition}\label{gal-duality}
	Let $\mathcal{P}_n$ be the set of prime numbers that split completely in $F_n(\ze{4})$. Let $2\le k\le n+1$. Then
	\[ \left(\frac{C_{F_n}}{\pm O_{F_n}^{\times 2^k}\cap C_{F_n}} \right)^{\perp}=\left\langle f_r\;:\;r\in\mathcal{P}_n \right\rangle_{R_k}.  \]
\end{proposition}

\begin{proof}
	By \cite[Theorem $6.3.2$]{gras}, the kernel of the natural map
	\[ \iota:O_{F_n}^{\times}/ O_{F_n}^{\times 2^k}\longrightarrow F_n(\ze{2^k})^\times/F_n(\ze{2^k})^{\times2^k} \]
	is cyclic of order $2$. Since $2\le k\le n+1$, we have $F_n(\ze{2^k})=F_n(\ze{4})=\Q(\ze{2^{n+2}},\sqrt{f})$. Moreover, $-1$ is a $2^k$-power in $F_n(\ze{4})$ but not in $F_n$ and hence $\ker\iota=\{\pm1\}$. So, we can identify $C_{F_n}/\pm O_{F_n}^{\times 2^k}\cap C_{F_n}$ with a subgroup of $F_n(\ze{4})^\times/F_n(\ze{4})^{\times2^k}$. Then, Kummer theory gives us an isomorphism
	\[ \gal\left(F_n(\ze{4})\left(\sqrt[2^k]{C_{F_n}}\right)|F_n(\ze{4})\right)\longrightarrow \HOM_{\Z}\left(\frac{C_{F_n}}{\pm O_{F_n}^{\times 2^k}\cap C_{F_n}},\frac{\Z}{2^k\Z}\right). \]
	This isomorphism depends on the choice of a $2^k$-primitive root of unity.
	On the other hand, since $R_k$ is finite Gorenstein we have the isomorphism described in Remark \ref{Gore-isomorphism}:
	\[ \HOM_{R_k}\left(\frac{C_{F_n}}{\pm O_{F_n}^{\times 2^k}\cap C_{F_n}},R_k\right)\longrightarrow \HOM_{\Z}\left(\frac{C_{F_n}}{\pm O_{F_n}^{\times 2^k}\cap C_{F_n}},\frac{\Z}{2^k\Z}\right). \]
	Choose a prime ideal $\mathfrak{R}\subset O_{F_n(\ze{4})}$ that divides $r\in\mathcal{P}_n$. It is easy to check that the image under the Kummer isomorphism of the Frobenius element associated to $\mathfrak{R}$ and the image of $f_{\mathfrak{R}}$ under the isomorphism given by the Gorenstein ring properties coincide.
	
	Recall that the Chebotar\"ev density theorem states that if we choose an automorphism $g\in\gal\left(F_n(\ze{4})\left(\sqrt[2^k]{C_{F_n}}\right)|F_n(\ze{4})\right)$, then there are infinitely many primes that split completely in $F_n(\ze{4})$ such that the Frobenius elements associated to those primes are equal to $g$. Therefore we get
	\[ \left(\frac{C_{F_n}}{\pm O_{F_n}^{\times 2^k}\cap C_{F_n}} \right)^{\perp}=\left\lbrace f_{\mathfrak{R}}\;:\; \mathfrak{R}|r \,\mathrm{and}\,r\in\mathcal{P}_n \right\rbrace = \left\langle f_r\;:\;r\in\mathcal{P}_n \right\rangle_{R_k}. \]
\end{proof}

Consider the exact sequence of $R_k$-modules presented in \ref{torsion}. We take the dual one:
\[ 0\longrightarrow\left(\frac{C_{F_n}}{\pm O_{F_n}^{\times 2^k}\cap C_{F_n}} \right)^{\perp}\longrightarrow\left(\frac{C_{F_n}}{\pm C_{F_n}^{ 2^k}} \right)^{\perp}\longrightarrow A_n[2^k]^{\perp}\longrightarrow0. \] 
We can use this sequence to compute the dual module of $A_n[2^k]$. But, before doing that we want to get a more manageable sequence.

\section{A short exact sequence suitable for computations}

In this section we find short exact sequences that we will use for computations. In particular, we describe the dual module of $A_n$ or $A_n/A_0$ as a quotient of a cyclic $\Z[G_n]$-module. We distinguish two cases according to the behaviour of the prime $2$ in the real quadratic field $F=\Q(\sqrt{f})$. The case in which the prime $2$ is inert or ramified in $F$ and the case in which the prime $2$ splits in $F$.

We consider now the  case in which the prime $2$ is ramified or inert in the real quadratic field  $F=\Q(\sqrt{f})$.

\begin{proposition}\label{nice-An-sequence1}
	Let $f\not\equiv1\bmod8$ be a positive odd squarefree integer. Let $F=\Q(\sqrt{f})$. Assume that the class number of $F$ is even.
	Denote by $X_{F_n}$ the $\Z[G_n]$-submodule of $C_{F_n}$ generated by $-1$ and $\eta_n$. Let $C_{\Q_n}$ be the $\Z[G_n]$-submodule of $C_{F_n}$ generated by $-1$ and $\beta_n$.
	For a fixed positive integer $k$ we denote by $Y_n$ the cokernel of the morphism
	\[ \frac{C_{\Q_n}}{\pm C_{\Q_n}^{2^k}}\longrightarrow
	\frac{C_{F_n}}{\pm O_{F_n}^{\times 2^k}\cap C_{F_n}}  \]
	induced by the inclusion $C_{\Q_n}\subset C_{F_n}$. Then, we have an exact sequence of $\Z[G_n]$-modules
	\[ 0\longrightarrow  A_n[2^k]\longrightarrow \frac{X_{F_n}}{\pm X_{F_n}^{2^k}}\longrightarrow Y_n\longrightarrow0.  \]				 
\end{proposition}

\begin{proof}
	Consider the natural commutative diagram with exact rows:
	\[ \xymatrix{0\ar[r]& 0\ar[r]\ar[d]&\frac{C_{\Q_n}}{\pm C_{\Q_n}^{2^k}}\ar[r]\ar[d]&\frac{C_{\Q_n}}{\pm C_{\Q_n}^{2^k}}\ar[r]\ar[d]&0\\
		0\ar[r]&  A_n[2^k]\ar[r]&\frac{C_{F_n}}{\pm C_{F_n}^{2^k}}\ar[r]&\frac{C_{F_n}}{\pm O_{F_n}^{\times2^k}\cap C_{F_n}}\ar[r]&0} \]
	where the bottom exact sequence is obtained as in \ref{torsion}.
	By Lemma \ref{Qn-injective} the rightmost map is injective.
	Moreover, 
	\[ \frac{C_{F_n}}{\pm1}=\frac{C_{\Q_n}}{\pm1}\oplus \frac{X_{F_n}}{\pm1} \]
	as a consequence of Theorem \ref{CFn-structure}. Then, applying the snake lemma we get
	\[ 0\longrightarrow  A_n[2^k]\longrightarrow \frac{X_{F_n}}{\pm X_{F_n}^{2^k}}\longrightarrow Y_n\longrightarrow0. \]
\end{proof}

Notice that the modules appearing in the previous proposition are modules over the ring $R_k=\Z/2^k\Z[G_n]$ and it is a finite Gorestein ring. Thus, taking the dual modules we get the following corollary.
\begin{corollary} \label{dual-ramified-inert}
	Let $f\not\equiv1\bmod8$ be an odd positive squarefree integer. Assume that the class number of $F=\Q(\sqrt{f})$ is even. Let $R_k=\Z/2^k\Z[G_n]$ for some positive integer $k$. Then
	\[ A_n[2^k]^{\perp}\cong \frac{ R_k}{\left\lbrace \psi(\eta_n)\;:\;\psi\in Y_n^\perp\right\rbrace}.  \]
	Moreover 
	\[Y_n^{\perp}=\left\lbrace \psi\in\left(\frac{C_{F_n}}{\pm O_{F_n}^{\times 2^k}\cap C_{F_n}}\right)^{\perp}\;:\; \psi(\beta_n)=0\right\rbrace .\]
\end{corollary}
\begin{proof}
	From the definition of $Y_n$ and the fact that $(-)^\perp$ is an exact functor on the category of finite $R_k$-modules we deduce that the sequence
	\[ 0\longrightarrow Y_n^\perp \longrightarrow \left( \frac{C_{F_n}}{\pm O_{F_n}^{\times 2^k}\cap C_{F_n}}\right)^{\perp} \longrightarrow \left( \frac{C_{\Q_n}}{\pm C_{\Q_n}^{\times 2^k}}\right)^{\perp} \]
	is exact. Since $C_{\Q_n}/\pm1$ is generated by $\beta_n$ we get that 
	\[Y_n^{\perp}=\left\lbrace \psi\in\left(\frac{C_{F_n}}{\pm O_{F_n}^{\times 2^k}\cap C_{F_n}}\right)^{\perp}\;:\; \psi(\beta_n)=0\right\rbrace .\]
	
	By Proposition \ref{nice-An-sequence1} and by the exactness of the functor $(-)^\perp$ we get an exact sequence of $R_k$-modules
	\[ 0\longrightarrow Y_n^\perp\longrightarrow \left( \frac{X_{F_n}}{\pm X_{F_n}^{2^k}}\right)^{\perp}\longrightarrow A_n[2^k]^\perp \longrightarrow 0.  \]	
	Due to Theorem \ref{CFn-structure} we have $X_{F_n}/\pm X_{F_n}^{2^k} \cong R_k$. Therefore we can identify $\left( X_{F_n}/\pm X_{F_n}^{2^k}\right) ^{\perp}$ with $ R_k$ via the isomorphism that maps $\psi$ to $\psi(\eta_n)$.  
	Via this identification we get the exact sequence
	\[ 0\longrightarrow\{\psi(\eta_n)\;:\;\psi\in Y_n^{\perp}\}\longrightarrow R_k\longrightarrow A_n[2^k]^{\perp}\longrightarrow 0. \]
\end{proof}

In the rest of this section we prove a result similar to Corollary \ref{dual-ramified-inert} in the case in which $2$ splits in the real quadratic field $F$.

\begin{proposition}\label{nice-An-sequence2}
	Let $f\ge3$ be a squarefree integer such that $f\equiv1\bmod8$. Denote by $X_{F_n}$ the $\Z[G_n]$-submodule of $C_{F_n}$ generated by $-1$ and $\eta_n$. Let $C_{\Q_n}$ be the $\Z[G_n]$-submodule of $C_{F_n}$ generated by $-1$ and $\beta_n$.  For a fixed positive integer $k$, we denote by $Y_n$ the cokernel of the natural morphism 
	\[ \frac{C_{\Q_n}}{\pm C_{\Q_n}^{2^k}}\oplus\frac{C_{F_0}}{\pm C_{F_0}^{2^k}}\longrightarrow \frac{C_{F_n}}{\pm O_{F_n}^{\times2^k}\cap C_{F_n}} \]
	induced by the inclusions $C_{\Q_n}\subset C_{F_n}$ and $C_{F_0}\subset C_{F_n}$. Then, we have an exact sequence of $\Z[G_n]$-modules
	\[ 0\longrightarrow \frac{A_n[2^k]}{A_0[2^k]}\longrightarrow\frac{X_{F_n}}{\pm X_{F_n}^{2^k}}\longrightarrow Y_n\longrightarrow0. \]
\end{proposition}
\begin{proof}
	We describe the $2^k$-torsion of $A_n$ with an exact sequence as in \ref{torsion}. In other words, we have
	\[0\longrightarrow  A_n[2^k]\longrightarrow \frac{C_{F_n}}{\pm C_{F_n}^{2^k}}\longrightarrow \frac{C_{F_n}}{\pm O_{F_n}^{\times2^k}\cap C_{F_n}}\longrightarrow0.\]
	
	In particular, the map
	\[ \frac{C_{\Q_n}}{\pm C_{\Q_n}^{2^k}}\oplus\frac{C_{F_0}}{\pm C_{F_0}^{2^k}}\longrightarrow \frac{C_{\Q_n}}{\pm C_{\Q_n}^{2^k}}\oplus\frac{C_{F_0}}{\pm O_{F_0}^{\times2^k}\cap C_{F_0}} \]
	induced by the identity on the first component and by the natural quotient on the second component has kernel isomorphic to $A_0[2^k]$.
	
	We consider the natural commutative diagram whose rows are exact:
	\[ \xymatrix{0\ar[r]&A_0[2^k]\ar[r]\ar[d]&\frac{C_{\Q_n}}{\pm C_{\Q_n}^{2^k}}\oplus \frac{C_{F_0}}{\pm C_{F_0}^{2^k}}\ar[r]\ar[d]&\frac{C_{\Q_n}}{\pm C_{\Q_n}^{2^k}}\oplus\frac{C_{F_0}}{\pm O_{F_0}^{\times 2^k}\cap C_{F_0}}\ar[r]\ar[d]&0\\
		0\ar[r]& A_n[2^k]\ar[r]&\frac{C_{F_n}}{\pm C_{F_n}^{2^k}}\ar[r]&\frac{C_{F_n}}{\pm O_{F_n}^{\times2^k}\cap C_{F_n}}\ar[r]&0.} \]
	
	We shown in Lemma \ref{Fn-Qn-injectivity} that the rightmost vertical map is injective. Furthermore, we have an isomorphism 
	\[ \frac{C_{F_n}}{\pm1}\cong\frac{X_{F_n}}{\pm1}\oplus\frac{C_{\Q_n}}{\pm1}\oplus\frac{C_{F_0}}{\pm1} \]
	as observed in Theorem \ref{CFn-structure}. Therefore, applying the snake lemma  we get the exact sequence
	\[ 0\longrightarrow \frac{A_n[2^k]}{A_0[2^k]}\longrightarrow\frac{X_{F_n}}{\pm X_{F_n}^{2^k}}\longrightarrow Y_n\longrightarrow0. \]
\end{proof}

\begin{corollary}\label{dual-split}
	Let $f\ge3$ be a squarefree integer such that $f\equiv1\bmod8$. Denote by $R_k$ the ring $\Z/2^k\Z[G_n]$ for some positive integer $k$. Then
	\[ \left( \frac{A_n[2^k]}{A_0[2^k]}\right) ^{\perp}\cong \frac{ I_k}{\left\lbrace \psi(\eta_n)\;:\;\psi\in Y_n^\perp\right\rbrace}\]
	where $I_k\subset R_k$ is the augmentation ideal. Moreover 
	\[Y_n^{\perp}=\left\lbrace \psi\in\left(\frac{C_{F_n}}{\pm O_{F_n}^{\times 2^k}\cap C_{F_n}} \right)^{\perp} \;:\;\psi(\beta_n)=\psi(\delta_f')=0 \right\rbrace.\]	
\end{corollary}
\begin{proof}
	The proof is essentially the same as the one of Corollary \ref{dual-ramified-inert}.
	By Proposition \ref{nice-An-sequence2} the module $Y_n^\perp$ fits in the exact sequences
	\begin{align*}
		0\longrightarrow Y_n^\perp \longrightarrow \left( \frac{X_{F_n}}{\pm X_{F_n}^{2^k}}\right)^{\perp} \longrightarrow \left( \frac{A_n[2^k]}{A_0[2^k]}\right) ^{\perp} \longrightarrow 0, \\
		0\longrightarrow Y_n^\perp \longrightarrow \left( \frac{C_{F_n}}{\pm O_{F_n}^{\times 2^k}\cap C_{F_n}}\right)^{\perp} \longrightarrow \left( \frac{C_{\Q_n}}{\pm C_{\Q_n}^{2^k}}\oplus\frac{C_{F_0}}{\pm C_{F_0}^{2^k}}\right)^{\perp} .		
	\end{align*}
	From the second one we deduce that $\psi\in\left(C_{F_n}/\pm O_{F_n}^{\times 2^k}\cap C_{F_n} \right)^{\perp}$ lies in $Y_n^\perp$ if and only if it is trivial when evaluated both in $\beta_n$ and $\delta_f'$. 
	
	Since $X_{F_n}/\pm X_{F_n}^{2^k}\cong R_k/(N_n)$ we have an isomorphism  $(X_{F_n}\pm X_{F_n}^{2^k})^\perp\cong I_k$ defined by mapping $\psi$ to $\psi(\eta_n)$. Under this identification $Y_n^\perp$ is mapped to
	\[ \left\lbrace \psi(\eta_n)\;:\;\psi\in Y_n^\perp\right\rbrace\subset I_k.\] 
\end{proof}

\section{The algorithm}\label{algo}

We explain in this section how we implement the above ideas on a computer in order to partially understand the module structure of $A_n$. In this section we fix a level $n$ and the integer $k=n+1$. Let $\pi$ be the the isomorphism
\begin{equation}\label{pi-iso} \tag*{($4.1$)}
	 \pi: R_{n+1}=\frac{\Z}{2^{n+1}\Z}[G_n]\longrightarrow\frac{\Z_2[X]}{(2^{n+1},X^{2^n}-1)}
\end{equation}
that maps $\gamma$ to $X$.
Via this isomorphism we translate the results of Corollaries \ref{dual-ramified-inert} and \ref{dual-split} in terms of polynomials so that a computer can easily handle them. For any $u\in C_{F_n}$ and $r\in\mathcal{P}_n$ we define the polynomial 
\[ f_r^u(X)=\sum_{i=0}^{2^n-1}\log_r(u^{\gamma^{-i}})X^i\in \frac{\Z_2[X]}{(2^{n+1},X^{2^n}-1)}. \]
Note that $\pi(f_r(u))=f_r^u(X)$.

Assume $f$ to be a positive odd squarefree integer such that $f\not\equiv1\bmod8$ and the class number of $\Q(\sqrt{f})$ is even. We want to find a collection of elements inside $Y_n^\perp$. 
We compute primes $r_1,r_2,\dots$ that are congruent to $1$ modulo $2^{n+2}f$. This implies that those primes split completely in $\Q(\ze{2^{n+2}f})$ and consequently in $F_n(\ze{4})$. Since $\beta_n\in C_{F_n}/\pm1$ is annihilated by $N_n$, the elements $f_{r_i}(\beta_n)$ lie in the augmentation ideal of $R_{n+1}$. In particular, $\gamma-1$ divides $f_{r_i}(\beta_n)$. For any pair $(r_i,r_j)$ we define
\[ g_{i,j}=\frac{f_{r_i}(\beta_n)}{\gamma-1}f_{r_j}-\frac{f_{r_j}(\beta_n)}{\gamma-1}f_{r_i}. \]
This is an element of $Y_n^{\perp}$. Indeed, $g_{i,j}(\beta_n)=0$. Let $J^{(1)}$ be the zero ideal of $R_{n+1}$. For any $i\ge2$, we compute the $R_{n+1}$-ideals
\[ J^{(i)}=J^{(i-1)}+(g_{i,j}(\eta_n)\;:\; 1\le j < i). \] 
Since $R_{n+1}$ is a finite ring, the ideals $J^{(i)}$ stabilize in finitely many steps to an ideal $J_n$. By Corollary \ref{dual-ramified-inert} there is a surjective homomorphism
\begin{equation}\label{surj-mor1} \tag*{($4.2$)}
\frac{R_{n+1}}{J_n}\longrightarrow A_n[2^{n+1}]^{\perp}.
\end{equation}

On a computer this strategy looks like the following: for a prime $r_i\equiv1\bmod 2^{n+2}f$ we compute the polynomials $f_{r_i}^{\beta_n}(X)$ and $f_{r_i}^{\eta_n}(X)$. Then, for any $j<i$ we have
\[ \pi(g_{i,j}(\eta_n))=\frac{f_{r_i}^{\beta_n}(X)}{X-1}f_{r_j}^{\eta_n}(X)-\frac{f_{r_j}^{\beta_n}(X)}{X-1}f_{r_i}^{\eta_n}(X). \]
If $M$ is an integer sufficiently large, then the morphism $\pi$ defined in \ref{pi-iso} induces an isomorphism
\[
\frac{R_{n+1}}{J_n}\overset{\pi}{\longrightarrow}\frac{\Z_2[X]}{\left(2^{n+1}, X^{2^{n}}-1,\pi(g_{i,j}(\eta_n))\;:\;1\le j<i\le M\right)}. 
\]

We compute ideals $J_n\subset R_{n+1}$ for $n=1,2,\dots$ until the hypotheses of the following Lemma are satisfied.

\begin{lemma}\label{end-algorithm}
	Let $f\not\equiv1\bmod8$ be an odd squarefree integer. Assume that the class number of $F=\Q(\sqrt{f})$ is even. Let $2^{m_0}$ be the cardinality of $A_0$. Assume that there exists an integer $m$ such that $2^m\in J_m$. Then, we have a surjective morphism of $\Z[G_m]$-modules
	\[ \frac{R_{m+1}}{J_m}\longrightarrow A_m^\perp. \] 
	Moreover, assume that one of the following properties holds.
	\begin{enumerate}[label=(\alph*)]
		\item\label{1.1-cond} The cardinality of $R_{m+1}/J_m$ is strictly less than $2^{m+m_0}$,
		\item\label{1.2-cond} the element $N_{m-1}^{(m)}\in R_{m+1}$ lies in the ideal $J_m$.
	\end{enumerate}
	Then, the natural morphisms $A_{m-1}\longrightarrow A_n$ are isomorphisms for any $n\ge m-1$. In particular, Greenberg's conjecture holds for $\Q(\sqrt{f})$.
	
\end{lemma} 
\begin{proof}
	We have the surjective  morphisms of $\Z[G_m]$-modules described in \ref{surj-mor1}
	\[ \frac{R_{m+1}}{J_m}\longrightarrow A_m[2^{m+1}]^\perp. \]
	Since $2^{m}$ lies in $J_{m}$ we get that $2^{m}$ annihilates $A_{m}[2^{m+1}]^{\perp}$. By Lemma \ref{dual-kill} the element $2^m$ annihilates $A_m[2^{m+1}]$ and Nakayama's lemma implies that $A_{m}[2^{m+1}]=A_{m}$. So, we are getting an upper bound for the cardinality of the full module $A_{m}$. 
	
	Assume that \ref{1.1-cond} holds. Then, we have
	\[ 2^{m_0}=\#A_0\le\#A_1\le\dots\le\#A_m<2^{m+m_0}. \]
	Therefore, there exists a level $m'<m$ so that $\#A_{m'}=\#A_{m'+1}$. Proposition \ref{An-stabilizatio} implies that $A_{m'}\cong A_n$ for any $n\ge m'$.
	
	Assume that \ref{1.2-cond} holds. Then $N_{m-1}^{(m)}$ annihilates $A_m^\perp$ and its submodule $(A_m/A_0)^\perp$. So, the element $N_{m-1}^{(m)}$ annihilates $A_m/A_0$ too. By Lemma \ref{relation-An}, we conclude that
	\[ \frac{A_{m-1}}{A_0}\cong\frac{A_m}{A_0}\left[N_{m-1}^{(m)}\right]=\frac{A_m}{A_0}. \]
	We get that $\#A_{m-1}=\#A_m$ since the natural morphisms $A_0\longrightarrow A_n$ are injective for any positive $n$. Finally, Proposition \ref{An-stabilizatio} implies that $A_{m-1}\cong A_n$ for any $n\ge m-1$. 
\end{proof}

Assume now that $f\ge3$ is a positive squarefree integer and $f\equiv1\bmod8$. We find several primes $r_1,r_2,\dots $ that are congruent to $1\bmod 2^{n+2}f$.  Since $\delta_f'$ is annihilated by $\gamma-1$, the elements $f_{r_i}(\delta_f')\in R_{n+1}$ are divisible by $N_n$. Fix $i\ge2$. Let $H_i$ be the subset of $(C_{F_n}/\pm O_{F_n}^{\time2^{n+1}}\cap C_{F_n})^{\perp}$ whose elements are
\[ \frac{f_{r_j}(\delta_f')}{2^{s(j,l)}N_n}f_{r_l}-\frac{f_{r_l}(\delta_f')}{2^{s(j,l)}N_n}f_{r_j} \] 
where $1\le j<l\le i$ and $s(j,l)\le n+1$ is the largest integer such that both $f_{r_j}(\delta_f')$ and $f_{r_l}(\delta_f')$ are congruent to $0$ modulo $2^{s(j,l)}$. The elements of $H_i$ are trivial when evaluated on $\delta_f'$. For any pair of elements $h_j,h_l\in H_i$ we define
\[ g_{h_j,h_l}=\frac{h_j(\beta_n)}{\gamma-1}h_l-\frac{h_l(\beta_n)}{\gamma-1}h_j. \]
Note that $g_{h_j,h_l}$ is trivial when evaluated on both $\delta_f'$ and $\beta_n$. Therefore, the functions $g_{h_j,h_l}$ lie in  $Y_n^\perp$. Let $J^{(i)}$ be the $R_{n+1}$-ideal 
\[ J^{(i)}=\left( g_{h_j,h_l}(\eta_n)\;:\;h_j,h_l\in H_i\right).  \]

It is easy to check that we have inclusions $J^{(i)}\subset J^{(i+1)}$. So the ideals $J^{(i)}$ stabilize in finitely many steps to an ideal $J_n$. Let $I_{n+1}\subset R_{n+1}$ be the augmentation ideal. By Corollary \ref{dual-split} there is surjective morphism
\begin{equation}\label{surj-mor2} \tag*{($4.3$)}
\frac{I_{n+1}}{J_n}\longrightarrow\left(\frac{A_n[2^{n+1}]}{A_0[2^{n+1}]}\right)^\perp.
\end{equation}
Furthermore, if $M$ is a sufficiently large integer then the morphism $\pi$ defined in \ref{pi-iso} induces an isomorphism
\[ \frac{I_{n+1}}{J_n}\overset{\pi}{\longrightarrow}\frac{(X-1)\Z_2[X]}{(2^{n+1}(X-1),X^{2^n}-1, \pi(g_{h_j,h_l}(\eta_n))\;:\;h_j,h_l\in H_M)}. \] 
Then, dividing by $X-1$ we get an isomorphism
\begin{equation}\label{split-quotient} \tag*{($4.4$)}
\frac{I_{n+1}}{J_n}{\longrightarrow}\frac{\Z_2[X]}{\left(2^{n+1},\frac{X^{2^n}-1}{X-1}, \frac{\pi(g_{h_j,h_l}(\eta_n))}{X-1}\;:\;h_j,h_l\in H_M\right)}
\end{equation} 
The polynomials $\pi(g_{h_j,h_l}(\eta_n))$ can be computed using combinations of $f_{r_i}^{\eta_n}(X)$, $f_{r_i}^{\beta_n}(X)$ and $f_{r_i}^{\delta_f'}(X)$.

For $n=1,2,\dots$ we compute the ideal $J_n$ until the hypotheses of the next Lemma are satisfied.

\begin{lemma}\label{end-algorithm2}
	Let $f\ge3$ be a squarefree integer such that $f\equiv1\bmod8$. Let $2^{m_0}$ be the cardinality of $A_0$. Assume that there exists an integer $m\ge m_0$ such that $2^{m-m_0}\in J_m$. Then, we have a surjective morphism of $\Z[G_m]$-modules
	\[ \frac{I_{m+1}}{J_m}\longrightarrow \left(\frac{A_m}{A_0}\right)^\perp. \] 
	Moreover, assume that one of the following properties holds.
	\begin{enumerate}[label=(\alph*)]
		\item\label{2.1-cond} The cardinality of $I_{m+1}/J_m$ is strictly less than $2^{m}$,
		\item\label{2.2-cond} the element $N_{m-1}^{(m)}\in R_{m+1}$ annihilates $I_{m+1}/J_m$.
	\end{enumerate}
	Then, the natural morphisms $A_{m-1}\longrightarrow A_n$ are isomorphisms for any $n\ge m-1$. In particular, Greenberg's conjecture holds for the field $\Q(\sqrt{f})$.
	
\end{lemma} 
\begin{proof}
	We have the surjective morphism of $\Z[G_m]$-modules described in \ref{surj-mor2}
	\[ \frac{I_{m+1}}{J_m}\longrightarrow\left( \frac{A_m[2^{m+1}]}{A_0}\right)^\perp. \]
	Since $2^{m-m_0}$ lies in $J_{m}$ we have that $2^{m-m_0}$ annihilates $(A_{m}[2^{m+1}]/A_0)^{\perp}$. By Lemma \ref{dual-kill} the element $2^{m-m_0}$ annihilates $A_m[2^{m+1}]/A_0$. It follows that $2^m$ annihilates $A_m[2^{m+1}]$ and Nakayama's lemma implies that $A_{m}[2^{m+1}]=A_{m}$.  
	
	Assume that \ref{1.1-cond} holds. Then, we have
	\[ 2^{m_0}=\#A_0\le\#A_1\le\dots\le\#A_m<2^{m+m_0}. \]
	Therefore, there exists a level $m'<m$ so that $\#A_{m'}=\#A_{m'+1}$. Proposition \ref{An-stabilizatio} implies that $A_{m'}\cong A_n$ for any $n\ge m'$.
	
	Assume that \ref{1.2-cond} holds. Then $N_{m-1}^{(m)}$ annihilates $(A_m/A_0)^\perp$. So, the element $N_{m-1}^{(m)}$ annihilates $A_m/A_0$ too. By Lemma \ref{relation-An}, we conclude that
	\[ \frac{A_{m-1}}{A_0}\cong\frac{A_m}{A_0}\left[N_{m-1}^{(m)}\right]=\frac{A_m}{A_0}. \]
	We get that $\#A_{m-1}=\#A_m$ since the natural morphisms $A_0\longrightarrow A_n$ are injective for any positive $n$. Finally, Proposition \ref{An-stabilizatio} implies that $A_{m-1}\cong A_n$ for any $n\ge m-1$. 
\end{proof}

We say that the algorithm terminates on an input $f$, if we find a level $m$ such that the hypothesis of Lemma \ref{end-algorithm} or Lemma \ref{end-algorithm2} accordingly are satisfied. We have proven the following.

\begin{theorem}\label{algorithm-terminates}
	If the algorithm terminates on an input $f$, then Greenberg's conjecture holds for $F=\Q(\sqrt{f})$.	
\end{theorem}

Happily, the algorithm succeeded in checking the conjecture for all the quadratic fields $\Q(\sqrt{f})$ such that $3\le f \le 10000$. Indeed, applying this algorithm we can not predict how many steps we need in order to find the level $m$ and the computational time for each step grows exponentially. Luckily, our computations never exceed level $12$ and the computer deals with them in reasonable time.
 
We address the problem of showing the converse of the above theorem in the next section while we conclude this section giving explicit formulae for the polynomials $f_{r}^{\eta_n}(X)$, $f_{r}^{\beta_n}(X)$ and $f_{r}^{\delta_f'}(X)$.
Assume $f\equiv3\bmod4$. Let $\chi_{-f}$ be the Dirichlet quadratic character for the field $\Q(\sqrt{-f})$. We have an isomorphism $\ker\chi_{-f}\cong \gal(\Q(\ze{f})|\Q(\sqrt{-f}))$. Notice that $F_{n+1}=\Q(\ze{2^{n+3}},\sqrt{-f})^+$, therefore we get
\[ \gal(\Q(\ze{2^{n+3}f})^+|F_{n+1})\cong \gal(\Q(\ze{2^{n+3}f})|\Q(\ze{2^{n+3}},\sqrt{-f}))\cong\ker\chi_{-f}. \]
We get that
\[ \eta_n=\prod_{a\in\ker\chi_{-f}}\left( \ze{4}(\ze{2^{n+3}}\ze{f}^a-\ze{2^{n+3}}^{-1}\ze{f}^{-a})\right).  \]
Moreover, a generator $\gamma$ of $G_n$ can be chosen as the inverse of the restriction to $F_n$ of $\phi\in\gal(\Q(\ze{2^{n+2}},\ze{f})|\Q(\ze{f}))$ defined by $\phi(\ze{2^{n+2}})=\ze{2^{n+2}}^3$. It follows that
\begin{equation}\label{freta} \tag*{($4.5$)} f_r^{\eta_n}(X)=\sum_{i=0}^{2^n-1}\log_r\left(\prod_{a\in\ker\chi_{-f}}\left(\ze{4}^{3^i}\left(\ze{2^{n+3}}^{3^i}\ze{f}^a-\ze{2^{n+3}}^{-3^i}\ze{f}^{-a}\right) \right)  \right)X^i.  
\end{equation}

On the other hand if $f\equiv1\bmod4$, let $\chi_f$ be the quadratic character for the field $\Q(\sqrt{f})$. As above we have that $F_{n+1}=\Q(\ze{2^{n+3}},\sqrt{f})^+$ and the group $\gal(\Q(\ze{2^{n+3}f})^+|F_{n+1})$ is isomorphic to $\ker\chi_f$. Therefore
\begin{equation}\label{fretabis} \tag*{($4.5'$)}
	 f_r^{\eta_n}(X)=\sum_{i=0}^{2^n-1}\log_r\left(\prod_{a\in\ker\chi_{f}}\left(\ze{4}^{3^i}\left(\ze{2^{n+3}}^{3^i}\ze{f}^a-\ze{2^{n+3}}^{-3^i}\ze{f}^{-a}\right) \right)  \right)X^i.  
\end{equation}
In both cases, we have
\[ f_r^{\beta_n}(X)=\sum_{i=0}^{2^n-1}\log_r\left(\ze{2^{n+2}}^{\frac{3^i-3^{i+1}}{2}}\frac{1-\ze{2^{n+2}}^{3^{i+1}}}{1-\ze{2^{n+2}}^{3^i}} \right)X^i. \]
Finally, if $f\equiv1\bmod 8$ we have
\[ f_r^{\delta_f'}(X)=f_r^{\delta_f^2}(X)=\log_r\left(\prod_{a\in\ker\chi_{f}}(1-\ze{f}^a)\right)\sum_{i=0}^{2^n-1}X^i\]
when $f$ is not prime. While if $f$ is prime we get
\[ f_r^{\delta_f'}(X)=f_r^{\delta_f}(X)=\log_r\left(\prod_{a\in\ker\chi_{f}/\pm1}\ze{f}^{a\frac{s-1}{2}}\frac{1-\ze{f}^{a}}{1-\ze{f}^{as}}\right)\sum_{i=0}^{2^n-1}X^i \]
where $s$ is an element in $(\Z/f\Z)^{\times}$ but not in $\ker\chi_f$. See Lemma \ref{square-properties} part \ref{prop-deltaf} for an explicit form of $\delta_f$.

\section{Further remarks on the algorithm and examples}\label{furtherrmk}

Notice that in the previous section we restrict our attention only to the auxiliary primes $r$ that are congruent to $1$ modulo $2^{n+2}f$ instead of the primes $r$ that splits in $F_n(\ze{4})$. In this section we discuss the reason why we do such assumption and what this implies.

We start explaining the reason. If we assume that $r\equiv1\bmod 2^{n+2}f$ we have that a primitive root of unity of order $2^{n+3}f$ lies in the finite field with $r^2$ elements, $\mathbb{F}_{r^2}$. Therefore, we can compute easily the product
\[ \prod_{a\in\ker\chi_{-f}}\left(\ze{4}^{3^i}\left(\ze{2^{n+3}}^{3^i}\ze{f}^a-\ze{2^{n+3}}^{-3^i}\ze{f}^{-a}\right) \right) \]
appearing in \ref{freta} or the one appearing in \ref{fretabis} by computing each single factor of the product in the finite field $\mathbb{F}_{r^2}$.

If we consider the weaker condition that $r$ splits in $F_n(\ze{4})$, then a primitive root of unity of order $2^{n+3}f$ lies in a finite extension of $\mathbb{F}_r$ whose degree divides $\varphi(f)$ where $\varphi$ denotes the Euleur totien function. So, computing the above product by computing each single factor in this finite extension of $\mathbb{F}_r$ is in general slower and in some cases we were not able to compute the product especially if $n$ grows and $f$ is large. 

In this way we are actually investigating a module that is bigger than $A_n$. Indeed, let $L_n=\Q(\ze{2^{n+2}f})^+$, then we define the $G_n$-module $E_n$ as
\[ E_n=\left\{u\in O_{L_n}^{\times}\;:\; u^{d}\in O_{F_n}^{\times}\t{ for some $d\in\Z_{>0}$}\right\} \]
and the quotient
\[ B_n=\frac{E_n}{C_{F_n}}\otimes_{\Z}\Z_2. \]
We can describe the $2$-torsion part of $B_n$ with an exact sequence similar to the one in \ref{torsion}. In other words, we have
\[ 0\longrightarrow B_n[2^k]\longrightarrow \frac{C_{F_n}}{\pm C_{F_n}^{2^k}}\longrightarrow\frac{C_{F_n}}{\pm O_{L_n}^{\times 2^k}\cap C_{F_n}}\longrightarrow 0 \]
and the dual sequence as $R_k$-modules
\[
 0\longrightarrow \left(\frac{C_{F_n}}{\pm O_{L_n}^{\times 2^k}\cap C_{F_n}}\right)^{\perp}\longrightarrow \left(\frac{C_{F_n}}{\pm C_{F_n}^{2^k}}\right)^{\perp}\longrightarrow B_n[2^k]^{\perp}\longrightarrow0.
\]

\begin{remark}\label{r1}
	Proceeding as for proving Proposition \ref{gal-duality} we can prove that
\[ \left(\frac{C_{F_n}}{\pm O_{L_n}^{\times 2^k}\cap C_{F_n}}\right)^{\perp}=\left\langle f_r\;:\; r\equiv1\bmod2^{n+2}f\right\rangle_{R_k}\subset \left(\frac{C_{F_n}}{\pm O_{F_n}^{\times 2^k}\cap C_{F_n}}\right)^{\perp}. \] 
Notice that in Proposition \ref{gal-duality} we consider primes $r$ that split in the field $F_n(\ze{4})$ while here we consider primes $r$ that splits in the field $L_n(\ze{4})=\Q(\ze{2^{n+2}f})$. So, the modules presented in the main tables are closely related to $B_n^\perp$ or $(B_n/B_0)^\perp$.
\end{remark}
However, although we restrict our attention to prime $r$ that are $1$ modulo $2^{n+2}f$, we still expect to be able to prove results about the stabilization of the modules $A_n$. Indeed, we have the following remarks:

\begin{remark}
	Consider the group morphism
	\[ E_n/O_{F_n}^{\times}\longrightarrow \HOM_{\Z}(\gal(L_n|F_n),\{\pm1\}) \]
	defined by $u\mapsto \phi_u$ where $\phi_u(\sigma)=u^{\sigma-1}$. It is well define. Indeed, $u^{\sigma-1}$ is a root of unity in $L_n$ which is a real field. Thus, we get that $\phi_u(\sigma)\in\{\pm1\}$. 
	
	Moreover, this morphism is injective and the group $\HOM_{\Z}(\gal(L_n|F_n),\{\pm1\})$ has exponent $2$. Therefore, we conclude that the exponent of $E_n/O_{F_n^{\times}}$ is at most $2$.
\end{remark}

\begin{remark}
	Recall that $\gal(L_n|F_n)$ can be identified with a subgroup of index $2$ of $(\Z/f\Z)^{\times}/\pm1$ if $n\ge1$. Therefore, we get
	\[ [B_n:A_n]=[E_n:O_{F_n}^{\times}]\le\#\HOM_{\Z}(\gal(L_n|F_n),\{\pm1\})\le 2^t \]
	where $t$ is the number of distinct primes that divide $f$. In particular, the modules $B_n$ stabilize in the $\Z_2$-extension if and only if the modules $A_n$ stabilize.
\end{remark}

We present now what are the main obstructions in proving the converse of Theorem \ref{algorithm-terminates}.
If we assume that Greenberg's conjecture holds for a quadratic field $F=\Q(\sqrt{f})$, we get that the modules $A_n$ and $B_n$ stabilize in the tower. Assume now that $2$ does not split in $F$. Let
\[ J_n^{\mathrm{full}}= \{f_r(\eta_n)\;:\; r\equiv1\bmod 2^{n+2}f \t{ and } f_r(\beta_n)=0\}.\]
Then, similarly to Corollary \ref{dual-ramified-inert}, it is easy to show that $B_n[2^{n+1}]^{\perp}\cong R_{n+1}/J_n^{\mathrm{full}}$.

The algorithm computes an $R_{n+1}$-ideal $J_n\subset J_n^{\mathrm{full}}$. If these two ideals are actually equal, and if $B_n$ stabilize in the tower then there exists a level $m$ such that the condition \ref{1.1-cond} of Lemma \ref{end-algorithm} is satisfied. So the algorithm would eventually terminate. But, since we are computing $J_n$ using only a finite number of functions $f_r$, we are not able to prove in general that $J_n=J_n^{\mathrm{full}}$. In order to solve this problem one could aim to further generalize the ideas of \cite[Section $4$]{kraftschoof}. However, we still cannot predict the number of steps the algorithm takes before finding the level $m$ and so the time taken to compute $m$ could be too long. 

On the other hand, if $2$ splits in $F$ an analogous situation occurs. Indeed, similarly to Corollary \ref{dual-split} we have that
\[ \left(\frac{B_n[2^{n+1}]}{B_0[2^{n+1}]}\right)^\perp\longleftarrow\frac{I_{n+1}}{\{f_r(\eta_n)\;:\; r\equiv1\bmod 2^{n+2}f \t{ and } f_r(\beta_n)=f_r(\delta_f')=0\}}. \]
is a surjective morphism with kernel of order at most $2$. This follows by the fact that the kernel of the natural morphism
\[ \frac{O_{L_0}^{\times}}{O_{L_0}^{\times 2^{n+1}}}\longrightarrow \frac{L_n^{\times}}{L_n^{\times 2^{n+1}}} \]
has order at most $2$ and we do not know if it is injective.

We conclude this section presenting two examples. The first one is an example of what usually happens. While the second one is an example where the group $A_n$ grows a little in the $\Z_2$-tower. 

\subsection*{Example} Let $f=949$. Then, $f\equiv5\mod8$ and the class number of $F=\Q(\sqrt{f})$ is $2$. Therefore, the cardinality of $A_0$ is $2$. For several $r\equiv1\bmod2^{n+2}f$ we compute the polynomials $f_r^{\eta_n}(X)$ and $f_r^{\beta_n}(X)$ in $\Z_2[X]/(2^{n+1},X^{2^n}-1)$. It is convenient to express such polynomials using the new parameter $T=X-1$.
Consider the level $n=1$.

\begin{longtable}{|c|c|c|}
	\hline
	$r$ & $f_r^{\eta_n}(T)$ & $f_r^{\beta_n}(T)/T$\\
	\hline
	\hline
	$22777$ & $T$ & $3$\\
	$45553$ & $3T+2$ & $3$\\
	$60737$ & $3T$ & $3$\\
	$68329$ & $2T + 2$ & $2$\\
	$136657$ & $3T$ & $1$\\
	$151841$ & $0$ & $2$\\
	\hline
\end{longtable}
The polynomials
\[ \frac{f_{r_i}^{\beta_n}(T)}{T}f_{r_j}^{\eta_n}(T)-\frac{f_{r_j}^{\beta_n}(T)}{T}f_{r_i}^{\eta_n}(T) \]
generate the ideal $(2)\subset\Z_2[T]/(4,(T-1)^2-1)$.
As in \ref{surj-mor1} we have a surjective homomorphism 
\[ \frac{\Z_2[T]}{(2,T^2)}\longrightarrow A_1[4]^{\perp}. \]
We deduce that $2$ annihilates $A_1[4]$. Hence $A_1[4]=A_1$.

Consider the level $n=2$. We compute the following polynomials

\begin{longtable}{|c|c|c|}
	\hline
	$r$ & $f_r^{\eta_n}(T)$ & $f_r^{\beta_n}(T)/T$\\
	\hline
	\hline
	$45553$ & $2T^3 + 5T + 2$ & $2T^2 + 2T +3$\\
	$60737$ & $7T^3 + 6T^2 + T + 4$ & $5T^2 +5T + 5$\\
	$136657$ & $4T^3 + 5T$ & $7T^2 +5T + 3$\\
	$151841$ & $T^3 + 7T^2 + 2T$ & $6T^2 +6T + 6$\\
	$182209$ & $6T^3 + 4T^2 + 2T + 4$ & $4T^2 + 7T$\\
	$273313$ & $6T^3 +2T^2 + 3T + 6$ & $5T^2 + 7T + 1$\\
	\hline
\end{longtable}
Note that the polynomials
\[ \frac{f_{r_i}^{\beta_n}(T)}{T}f_{r_j}^{\eta_n}(T)-\frac{f_{r_j}^{\beta_n}(T)}{T}f_{r_i}^{\eta_n}(T) \]
generate the ideal $(2,T^2)\subset\Z_2[T]/(8,(T-1)^4-1)$. So, we get a surjective morphism
\[ \frac{\Z_2[T]}{(2,T^2)}\longrightarrow A_2[8]^{\perp}. \]
We conclude that $A_2[8]=A_2$ and $A_2\le4$. The hypotheses of Lemma \ref{end-algorithm} are satisfied. We can deduce that $A_1\cong A_n$ for any $n\ge1$. 

Using Pari-GP, we see that the class group $\cl_{F_1}$ is cyclic of order $2$. Therefore, the $2$-parts of the class groups actually stabilize at level $0$. 

\subsection*{Example} Let $f=6817$. In this case, $f\equiv1\bmod8$ and the field $F=\Q(\sqrt{f})$ has class number $2$. For $n=1,2,\dots$ we compute ideals $J_n'$ of $\Z_2[X]$ such that we have a surjective morphism
\[ \frac{\Z_2[X]}{J_n'}\longrightarrow\left(\frac{A_n[2^{n+1}]}{A_0[2^{n+1}]}\right)^{\perp}. \]
The ideals $J_n'$ are the ones appearing in \ref{split-quotient}. They are computed using $15$ auxiliary primes $r$. We collect below these ideals using the parameter $T=X-1$ and the times taken for computing them. We see that For $n=7$ the hypotheses of Lemma \ref{end-algorithm2} are satisfied. Indeed, the polynomials $2^6$ and $\sum_{i=0}^{2^6-1}(T+1)^i$ lies in $J_7'$. Note that via the isomorphism described in \ref{split-quotient} the polynomial $\sum_{i=0}^{2^6-1}(T+1)^i$ is mapped to $N_6^{(7)}=\sum_{i=0}^{2^6-1}\gamma^i$. Therefore we get that for any $n\ge6$ the natural morphisms $A_6\longrightarrow A_n$ are isomorphisms. Notice that in this case condition \ref{2.2-cond} of Lemma \ref{end-algorithm2} permit us to stop our computations at level $7$. While condition \ref{2.1-cond} would have require to extend the computation up to level $11$.

\newpage
\begin{longtable}{|c|c|c|c|}
	\hline
	$n$ & $J_n'\subset\Z_2[T]$ & $[\Z_2[T]:J_n']$ & time\\
	\hline
	\hline
	$1$ & $(4, T + 2)$ & $4$ & $1.29s$\\
	$2$ & $(8, 4T, 2T^2, T^3 + 2T + 4)$ & $64$ & $2.42s$\\
	$3$ & $(8, 4T, 2T^2, T^4)$ & $128$ & $4.80s$\\
	$4$ & $(16, 4T, 2T^2, T^4)$ & $256$ & $9.58s$\\
	$5$ & $(32, 4T, 2T^2, T^4)$ & $512$ & $19.51s$\\
	$6$ & $(64, 4T, 2T^2, T^4 + 32)$ & $1024$ & $40.99s$\\
	$7$ & $(64, 4T, 2T^2, T^4 + 32)$ & $1024$ & $83.95s$\\	
	\hline
\end{longtable}

We compare the results with the ones we obtain using Pari-GP. We find out that the class groups $\cl_{F_1}$, $\cl_{F_2}$, $\cl_{F_3}$ are isomorphic as groups respectively to $\Z/4\Z\times\Z/2\Z$, $\Z/8\Z\times\Z/2\Z\times\Z/2\Z\times\Z/2\Z$ and $\Z/8\Z\times\Z/4\Z\times\Z/2\Z\times\Z/2\Z$. The Pari packages are faster in computing the class group at the first and the second level taking only $0.001s$ and $0.467s$, but it becomes slower from the third level. It takes $38.67s$ for computing $\cl_{F_3}$. We do not compute $\cl_{F_4}$ since the computational time exceeds one hour.

\section{The main tables}

We collect now the results we get. In all the cases we found that the modules $A_n$ stabilize. As a consequence we get that Greenberg's conjecture hold in all cases. The computations were done using SageMath $8.8$, the algorithm we use is available at \url{https://drive.google.com/drive/folders/1vNyOPfQadLd9yFk4qCOfDm24lJuPWlhQ?usp=sharing}.

\subsection*{Results for $f\equiv3\bmod4$} Let $3\le f<10000$ be a squarefree integer so that $f$ is congruent to $3$ modulo $4$. Assume that the class number of $F_0=\Q(\sqrt{f})$ is even. As presented in Section \ref{algo}, we compute an ideal $J\subset\Z_2[X]$ such that for any positive $n$ we have a surjective morphism
\[ \frac{\Z_2[X]}{J}\longrightarrow A_n^{\perp}.\]
Denote by $N$ the index $[\Z_2[X]:J]$.
Let $n_0$ be the least integer such that the polynomial $\sum_{i=0}^{2^{n_0}-1}X^i$ lies in $J$. By Lemma \ref{relation-An}, the natural morphisms $A_{n_0}\longrightarrow A_n$ are isomorphisms for all $n\ge n_0$.  

We collect the ideals $J$ in the following table. It is convenient to write the results with the parameter $T=X-1$.

\begingroup
\fontsize{8}{10}\selectfont

\begin{longtable}{|p{0.33\textwidth}|p{0.02\textwidth}|p{0.02\textwidth}|p{0.52\textwidth}|}
	
	\hline
	$J$ & $n_0$ & $N$ & $f$ \\ 
	\hline
	\hline 
	\endhead
	
	\hline
	\endfoot
	
	$ (2,T) $ & $ 1 $ & $ 2 $ & $15$, $35$, $39$, $55$, $87$, $91$, $95$, $111$, $115$, $143$, $155$, $159$, $183$, $203$, $215$, $235$, $247$, $259$, $295$, $299$, $303$, $319$, $327$, $335$, $355$, $371$, $395$, $403$, $407$, $415$, $427$, $447$, $471$, $515$, $519$, $535$, $543$, $551$, $559$, $583$, $591$, $611$, $635$, $655$, $667$, $671$, $687$, $695$, $703$, $707$, $755$, $763$, $767$, $807$, $815$, $831$, $835$, $851$, $871$, $879$, $895$, $899$, $923$, $951$, $955$, $995$, $1007$, $1027$, $1043$, $1047$, $1055$, $1079$, $1099$, $1111$, $1115$, $1119$, $1135$, $1147$, $1159$, $1167$, $1191$, $1195$, $1199$, $1211$, $1219$, $1247$, $1255$, $1263$, $1267$, $1315$, $1339$, $1355$, $1363$, $1379$, $1383$, $1391$, $1403$, $1415$, $1527$, $1535$, $1555$, $1591$, $1603$, $1623$, $1639$, $1643$, $1651$, $1655$, $1671$, $1703$, $1711$, $1727$, $1735$, $1739$, $1795$, $1807$, $1835$, $1839$, $1883$, $1891$, $1895$, $1903$, $1915$, $1919$, $1939$, $1943$, $1959$, $1963$, $1983$, $1991$, $2031$, $2051$, $2059$, $2071$, $2095$, $2103$, $2119$, $2127$, $2155$, $2167$, $2171$, $2183,\dots$
	\\ & & & $\dots2195$, $2199$, $2215$, $2219$, $2271$, $2279$, $2291$, $2315$, $2319$, $2323$, $2327$, $2335$, $2391$, $2395$, $2407$, $2435$, $2443$, $2455$, $2463$, $2479$, $2483$, $2487$, $2491$, $2495$, $2507$, $2515$, $2519$, $2559$, $2587$, $2611$, $2615$, $2623$, $2627$, $2631$, $2723$, $2735$, $2743$, $2779$, $2815$, $2823$, $2831$, $2855$, $2867$, $2899$, $2923$, $2935$, $2947$, $2951$, $2959$, $2983$, $2987$, $2991$, $2995$, $3035$, $3039$, $3047$, $3063$, $3071$, $3095$, $3103$, $3107$, $3127$, $3131$, $3155$, $3183$, $3207$, $3215$, $3223$, $3227$, $3235$, $3263$, $3279$, $3287$, $3295$, $3327$, $3351$, $3379$, $3415$, $3419$, $3427$, $3439$, $3455$, $3487$, $3523$, $3543$, $3551$, $3563$, $3595$, $3599$, $3611$, $3635$, $3639$, $3679$, $3683$, $3687$, $3695$, $3711$, $3715$, $3743$, $3755$, $3763$, $3787$, $3799$, $3811$, $3831$, $3839$, $3899$, $3903$, $3935$, $3959$, $3979$, $3991$, $4031$, $4043$, $4055$, $4087$, $4103$, $4115$, $4119$, $4135$, $4143$, $4163$, $4187$, $4195$, $4279$, $4287$, $4291$, $4295$, $4303$, $4315$, $4331$, $4343$, $4351$, $4359$, $4367$, $4379$, $4399$, $4415$, $4435$, $4479$, $4511$, $4531$, $4535$, $4555$, $4571$, $4595$, $4619$, $4627$, $4631$, $4647$, $4667$, $4687$, $4699$, $4727$, $4735$, $4739$, $4747$, $4771$, $4791$, $4819$, $4835$, $4839$, $4843$, $4847$, $4855$, $4863$, $4867$, $4907$, $4911$, $4915$, $4927$, $4955$, $4963$, $4979$, $5007$, $5063$, $5071$, $5079$, $5095$, $5111$, $5123$, $5127$, $5131$, $5143$, $5155$, $5191$, $5195$, $5199$, $5223$, $5255$, $5263$, $5267$, $5299$, $5315$, $5363$, $5367$, $5411$, $5435$, $5447$, $5455$, $5459$, $5515$, $5539$, $5567$, $5579$, $5583$, $5587$, $5599$, $5603$, $5611$, $5615$, $5631$, $5671$, $5703$, $5707$, $5747$, $5755$, $5759$, $5771$, $5799$, $5803$, $5815$, $5847$, $5855$, $5919$, $5935$, $5951$, $5959$, $5971$, $5991$, $6019$, $6023$, $6031$, $6071$, $6087$, $6107$, $6115$, $6119$, $6127$, $6139$, $6155$, $6159$, $6179$, $6187$, $6207$, $6227$, $6283$, $6295$, $6331$, $6371$, $6383$, $6395$, $6407$, $6415$, $6423$, $6431$, $6455$, $6467$, $6487$, $6515$, $6527$, $6535$, $6539$, $6583$, $6587$, $6595$, $6623$, $6631$, $6635$, $6639$, $6663$, $6711$, $6731$, $6739$, $6743$, $6751$, $6767$, $6799$, $6807$, $6835$, $6879$, $6927$, $6931$, $6943$, $6979$, $6995$, $6999$, $7003$, $7023$, $7067$, $7071$, $7087$, $7091$, $7099$, $7111$, $7115$, $7135$, $7143$, $7147$, $7167$, $7171$, $7183$, $7195$, $7235$, $7255$, $7271$, $7279$, $7291$, $7295$, $7303$, $7311$, $7319$, $7355$, $7363$, $7367$, $7379$, $7391$, $7415$, $7423$, $7427$, $7431$, $7435$, $7439$, $7447$, $7483$, $7495$, $7543$, $7555$, $7615$, $7627$, $7631$, $7647$, $7651$, $7655$, $7671$, $7711$, $7715$, $7739$, $7747$, $7763$, $7783$, $7787$, $7795$, $7799$, $7807$, $7819$, $7835$, $7855$, $7859$, $7863$, $7891$, $7895$, $7915$, $7979$, $7991$, $7999$, $8003$, $8027$, $8031$, $8035$, $8047$, $8063$, $8079$, $8095$, $8131$, $8135$, $8203$, $8207$, $8223$, $8247$, $8251$, $8267$, $8315$, $8327$, $8335$, $8339$, $8359$, $8367$, $8383$, $8391$, $8399$, $8411$, $8471$, $8479$, $8491$, $8495$, $8503$, $8507$, $8511$, $8567$, $8579$, $8583$, $8587$, $8603$, $8611$, $8615$, $8639$, $8659$, $8727$, $8735$, $8751$, $8759$, $8767$, $8791$, $8795$, $8843$, $8851$, $8871$, $8879$, $8903$, $8915$, $8935$, $8939$, $8947$, $8983$, $9019$, $9031$, $9047$, $9055$, $9083$, $9107$, $9111$, $9115$, $9119$, $9131$, $9155$, $9183$, $9211$, $9235$, $9259$, $9263$, $9287$, $9327$, $9335$, $9347$, $9355$, $9383$, $9395$, $9451$, $9487$, $9535$, $9543$, $9599$, $9607$, $9611$, $9647$, $9655$, $9659$, $9663$, $9667$, $9671$, $9683$, $9687$, $9731$, $9755$, $9759$, $9763$, $9827$, $9847$, $9895$, $9903$, $9935$, $9943$, $9983$, $9995$.\\ 
	$ (2,T^2) $ & $ 2 $ & $ 2^2 $ & $51$, $123$, $187$, $287$, $411$, $451$, $511$, $623$, $779$, $1203$, $1347$, $1563$, $1707$, $1819$, $1843$, $1923$, $1927$, $2047$, $2123$, $2191$, $2227$, $2263$, $2283$, $2363$, $2427$, $2571$, $2651$, $2747$, $2759$, $2771$, $2787$, $2859$, $2863$, $2911$, $2931$, $3151$, $3239$, $3431$, $3587$, $3647$, $3667$, $3859$, $4083$, $4207$, $4227$, $4247$, $4499$, $4579$, $4811$, $5027$, $5091$, $5327$, $5723$, $5891$, $5899$, $5999$, $6191$, $6243$, $6259$, $6439$, $6443$, $6459$, $6463$, $6499$, $6559$, $6847$, $6891$, $7123$, $7199$, $7251$, $7339$, $7343$, $7403$, $7519$, $7619$, $7679$, $7771$, $7827$, $7831$, $7899$, $8051$, $8159$, $8187$, $8227$, $8259$, $8331$, $8351$, $8371$, $8403$, $8459$, $8483$, $8531$, $8651$, $8691$, $8899$, $9123$, $9167$, $9179$, $9247$, $9299$, $9307$, $9427$, $9627$, $9703$, $9707$, $9727$, $9799$, $9899$, $9987$.\\ 
	$ (4,2T,T^2) $ & $ 2 $ & $ 2^3 $ & $195$, $219$, $231$, $291$, $555$, $579$, $715$, $731$, $759$, $903$, $915$, $939$, $943$, $959$, $979$, $1003$, $1015$, $1131$, $1139$, $1227$, $1235$, $1271$, $1295$, $1311$, $1371$, $1407$, $1411$, $1463$, $1495$, $1551$, $1595$, $1631$, $1635$, $1679$, $1767$, $1803$, $1855$, $1967$, $2015$, $2067$, $2247$, $2307$, $2343$, $2355$, $2419$, $2603$, $2607$, $2639$, $2679$, $2715$, $2755$, $2795$, $2967$, $3027$, $3199$, $3219$, $3335$, $3355$, $3367$, $3423$, $3435$, $3443$, $3515$, $3759$, $3815$, $3835$, $3891$, $3939$, $4155$, $4171$, $4183$, $4191$, $4223$, $4267$, $4319$, $4411$, $4427$, $4431$, $4467$, $4495$, $4503$, $4543$, $4615$, $4623$, $4763$, $4823$, $4891$, $4939$, $4983$, $5035$, $5207$, $5215$, $5235$, $5291$, $5307$, $5371$, $5403$, $5487$, $5555$, $5595$, $5627$, $5663$, $5699$, $5719$, $5727$, $5731$, $5735$, $5767$, $5811$, $5871$, $5883$, $5955$, $5963$, $5979$, $5995$, $6063$, $6231$, $6303$, $6315$, $6319$, $6391$, $6411$, $6479$, $6483$, $6695$, $6747$, $6787$, $6843$, $6895$, $6951$, $6955$, $7031$, $7131$, $7223$, $7231$, $7239$, $7287$, $7383$, $7419$, $7531$, $7579$, $7683$, $7903$, $7959$, $8115$, $8139$, $8195$, $8215$, $8255$, $8571$, $8607$, $8695$, $8711$, $8859$, $8891$, $8911$, $8943$, $9147$, $9195$, $9271$, $9303$, $9407$, $9447$, $9455$, $9483$, $9507$, $9515$, $9571$, $9695$, $9699$, $9715$, $9815$, $9915$, $9919$, $9951$.\\ 
	$ (2,T^3) $ & $ 2 $ & $ 2^3 $ & $119$, $339$, $391$, $679$, $771$, $1059$, $1207$, $1687$, $2147$, $2231$, $2839$, $3031$, $3383$, $4803$, $4859$, $4883$, $5287$, $5383$, $5543$, $6503$, $6667$, $6819$, $6887$, $7463$, $7571$, $7971$, $8279$, $8551$, $9223$, $9267$, $9411$, $9527$, $9959$.\\ 
	$ (4,2T,T^2 + 2) $ & $ 2 $ & $ 2^3 $ & $267$, $803$, $843$, $1691$, $3291$, $3579$, $3867$, $4299$, $4307$, $6059$, $7179$, $9771$.\\ 
	$ (4,2T,T^3) $ & $ 2 $ & $ 2^4 $ & $483$, $615$, $651$, $663$, $1095$, $1435$, $1455$, $1491$, $1599$, $1615$, $1659$, $1751$, $1771$, $2055$, $2139$, $2255$, $2387$, $2431$, $2555$, $2567$, $2667$, $2703$, $2847$, $2895$, $3115$, $3143$, $3471$, $3495$, $3507$, $3567$, $3619$, $3655$, $3731$, $3783$, $3895$, $4015$, $4439$, $4471$, $4551$, $4683$, $4711$, $5015$, $5151$, $5335$, $5343$, $5451$, $5467$, $5523$, $5691$, $5695$, $5863$, $6015$, $6251$, $6351$, $6495$, $6519$, $6531$, $6603$, $6643$, $6735$, $6839$, $6923$, $6935$, $7055$, $7107$, $7503$, $7527$, $7539$, $7599$, $7707$, $7743$, $7815$, $7843$, $7931$, $8007$, $8099$, $8323$, $8439$, $8535$, $8763$, $9255$, $9331$, $9367$, $9399$, $9443$, $9499$, $9579$, $9615$, $9635$, $9723$, $9863$, $9879$, $9911$.\\ 
	$ (8,2T + 4,T^2 + 4) $ & $ 2 $ & $ 2^4 $ & $399$, $1299$, $1443$, $1743$, $2035$, $2379$, $2751$, $2915$, $2919$, $3399$, $3403$, $3963$, $3999$, $4387$, $4767$, $5135$, $5219$, $5271$, $5947$, $6123$, $6447$, $6567$, $6671$, $7563$, $7719$, $7847$, $8555$, $8671$, $9023$, $9035$, $9143$, $9595$, $9807$, $9955$.\\ 
	$ (4,T^2 + 2) $ & $ 2 $ & $ 2^4 $ & $699$, $3091$, $3139$, $3827$, $5163$, $5739$, $6611$, $7811$.\\ 
	$ (4,2T^2,T^3 + 2T) $ & $ 2 $ & $ 2^5 $ & $255$, $935$, $987$, $1479$, $1887$, $3111$, $4123$, $4215$, $4795$, $5019$, $5423$, $5559$, $6919$, $8043$, $8155$, $8455$, $8499$, $8823$, $9051$, $9095$, $9215$, $9231$, $9379$, $9939$.\\ 
	$ (16,2T + 4,T^2 + 12) $ & $ 2 $ & $ 2^5 $ & $4251$, $7511$, $9823$.\\ 
	$ (8,2T + 4,T^3) $ & $ 2 $ & $ 2^5 $ & $1351$.\\ 
	$ (8,4T,2T^2 + 4,T^3 + 2T + 4) $ & $ 2 $ & $ 2^6 $ & $2827$, $3747$, $6707$, $9219$.\\ 
	$ (16,2T + 4,T^3 + 8) $ & $ 2 $ & $ 2^6 $ & $7063$.\\ 
	$ (8,2T^2 + 4T + 4,T^3 + 6T + 4) $ & $ 2 $ & $ 2^7 $ & $4695$, $6855$.\\ 
	$ (8,2T,T^2) $ & $ 3 $ & $ 2^4 $ & $323$, $435$, $455$, $723$, $795$, $1507$, $1515$, $1763$, $2955$, $3043$, $3055$, $3387$, $4035$, $4071$, $4147$, $4395$, $4755$, $5159$, $5395$, $5495$, $5511$, $5551$, $5619$, $6051$, $6055$, $6235$, $6335$, $6815$, $7015$, $7635$, $7851$, $7887$, $7939$, $8015$, $8239$, $8319$, $9003$, $9039$, $9139$, $9191$, $9415$, $9519$, $9591$, $9979$.\\ 
	$ (2,T^4) $ & $ 3 $ & $ 2^4 $ & $527$, $1343$, $1799$, $2471$, $3651$, $3707$, $4151$, $4659$, $5311$, $6167$, $6239$, $6387$, $6403$, $7327$, $7471$, $7751$, $8071$, $8143$, $8743$, $8927$.\\ 
	$ (4,2T,T^4) $ & $ 3 $ & $ 2^5 $ & $791$, $1011$, $1547$, $2359$, $3503$, $3603$, $4407$, $5083$, $6215$, $6307$, $6755$, $8407$, $8435$, $8755$.\\ 
	$ (16,2T + 12,T^2 + 12) $ & $ 3 $ & $ 2^5 $ & $1239$, $2019$, $4971$, $5795$, $5943$, $6771$, $7059$, $8083$, $8515$, $8635$, $8683$, $8799$.\\ 
	$ (8,4T,T^2 + 2T + 4) $ & $ 3 $ & $ 2^5 $ & $1419$, $3099$, $4587$, $4731$, $5259$, $6963$, $7467$, $7611$, $9339$.\\ 
	$ (8,4T,T^2 + 2T) $ & $ 3 $ & $ 2^5 $ & $1947$, $3531$, $3819$, $5379$, $5907$, $8643$, $8987$.\\ 
	$ (8,4T,T^2 + 6) $ & $ 3 $ & $ 2^5 $ & $2563$, $4443$, $5251$, $7387$.\\ 
	$ (8,2T,T^3) $ & $ 3 $ & $ 2^5 $ & $2807$, $4487$, $6103$, $9079$.\\ 
	$ (4,2T^2,T^3) $ & $ 3 $ & $ 2^5 $ & $3459$, $5331$, $9363$, $9651$.\\ 
	$ (2,T^5) $ & $ 3 $ & $ 2^5 $ & $4039$, $5667$, $8519$, $9071$.\\ 
	$ (4,2T,T^4 + 2) $ & $ 3 $ & $ 2^5 $ & $1779$, $2643$, $8023$.\\ 
	$ (8,4T,2T^2,T^3) $ & $ 3 $ & $ 2^6 $ & $1155$, $1995$, $3135$, $3615$, $3795$, $4199$, $4515$, $4895$, $5115$, $5187$, $5655$, $6135$, $6195$, $6355$, $7035$, $7095$, $7755$, $8103$, $8151$, $8547$, $9015$.\\ 
	$ (4,2T^2,T^4 + 2T) $ & $ 3 $ & $ 2^6 $ & $595$, $1695$, $1955$, $2159$, $3451$, $4403$, $6035$, $6339$, $7259$, $9831$.\\ 
	$ (4,2T,T^5) $ & $ 3 $ & $ 2^6 $ & $2635$, $3955$, $4371$, $4559$, $6715$, $6851$, $7347$, $7663$, $8655$, $8895$.\\ 
	$ (4,2T^2,T^4) $ & $ 3 $ & $ 2^6 $ & $2091$, $3723$, $4991$, $6987$, $7567$, $7667$.\\ 
	$ (8,2T + 4,T^4) $ & $ 3 $ & $ 2^6 $ & $799$, $4063$, $4607$.\\ 
	$ (8,4T,2T^2 + 4,T^3 + 2T) $ & $ 3 $ & $ 2^6 $ & $1335$, $8395$.\\ 
	$ (32,2T + 20,T^2 + 28) $ & $ 3 $ & $ 2^6 $ & $4807$, $7359$.\\ 
	$ (16,4T + 8,T^2 + 2T + 8) $ & $ 3 $ & $ 2^6 $ & $2211$.\\ 
	$ (8,T^2 + 6) $ & $ 3 $ & $ 2^6 $ & $3147$.\\ 
	$ (8,2T,T^4) $ & $ 3 $ & $ 2^6 $ & $3791$.\\ 
	$ (4,2T^2,T^4 + 2) $ & $ 3 $ & $ 2^6 $ & $6523$.\\ 
	$ (8,T^2 + 4T + 2) $ & $ 3 $ & $ 2^6 $ & $7323$.\\ 
	$ (16,2T + 12,T^3 + 8) $ & $ 3 $ & $ 2^6 $ & $9991$.\\ 
	$ (8,4T,2T^2,T^4) $ & $ 3 $ & $ 2^7 $ & $3255$, $4935$, $7995$, $8855$, $9435$, $9471$.\\ 
	$ (16,4T + 8,2T^2 + 8,T^3 + 8) $ & $ 3 $ & $ 2^7 $ & $7315$, $8815$, $8835$, $9735$, $9779$, $9867$.\\ 
	$ (8,2T^2 + 4,T^3 + 2T) $ & $ 3 $ & $ 2^7 $ & $1243$, $2163$, $4011$, $4179$.\\ 
	$ (4,2T^3,T^4 + 2T^2 + 2T) $ & $ 3 $ & $ 2^7 $ & $3395$, $6511$, $8827$.\\ 
	$ (4,2T^2,T^5 + 2T) $ & $ 3 $ & $ 2^7 $ & $3995$, $5055$.\\ 
	$ (16,4T + 8,2T^2,T^3 + 8) $ & $ 3 $ & $ 2^7 $ & $3883$.\\ 
	$ (8,4T,2T^2 + 4,T^4 + 4) $ & $ 3 $ & $ 2^7 $ & $4539$.\\ 
	$ (4,2T^3,T^4 + 2T) $ & $ 3 $ & $ 2^7 $ & $5295$.\\ 
	$ (4,2T^3,T^4 + 2T^2) $ & $ 3 $ & $ 2^7 $ & $6347$.\\ 
	$ (16,4T + 8,2T^2,T^3) $ & $ 3 $ & $ 2^7 $ & $7779$.\\ 
	$ (8,4T,2T^3,T^4) $ & $ 3 $ & $ 2^8 $ & $6279$, $7455$.\\ 
	$ (8,2T^2 + 4,T^4 + 4T + 4) $ & $ 3 $ & $ 2^8 $ & $8067$, $8119$.\\ 
	$ (4,T^4 + 2T^3 + 2T^2 + 2) $ & $ 3 $ & $ 2^8 $ & $2599$.\\ 
	$ (8,4T,2T^3,T^4 + 2T^2) $ & $ 3 $ & $ 2^8 $ & $3315$.\\ 
	$ (8,4T,2T^2,T^5) $ & $ 3 $ & $ 2^8 $ & $3927$.\\ 
	$ (4,2T^2,T^6) $ & $ 3 $ & $ 2^8 $ & $5763$.\\ 
	$ (16,4T + 8,2T^2 + 8,T^4) $ & $ 3 $ & $ 2^8 $ & $7395$.\\ 
	$ (4,2T^3,T^5 + 2T) $ & $ 3 $ & $ 2^8 $ & $8995$.\\ 
	$ (8,4T,2T^3,T^5 + 2T^2) $ & $ 3 $ & $ 2^9 $ & $6783$, $7735$.\\ 
	$ (4,2T^3,T^6 + 2T^2) $ & $ 3 $ & $ 2^9 $ & $9843$.\\ 
	$ (4,2T^4,T^6 + 2T^2) $ & $ 3 $ & $ 2^{10} $ & $4947$.\\ 
	$ (8,4T,2T^3,T^6) $ & $ 3 $ & $ 2^{10} $ & $8211$.\\ 
	$ (16,8T,4T^2 + 8,2T^3 + 4T + 8,T^4 + 4T + 4) $ & $ 3 $ & $ 2^{10} $ & $9691$.\\ 
	$ (16,2T + 8,T^2) $ & $ 4 $ & $ 2^5 $ & $2135$, $2235$, $3311$, $3535$, $4255$, $5359$, $8355$, $9159$.\\ 
	$ (16,2T,T^2) $ & $ 4 $ & $ 2^5 $ & $3983$, $4047$, $4355$, $8787$, $9795$.\\ 
	$ (32,2T + 12,T^2 + 28) $ & $ 4 $ & $ 2^6 $ & $6683$, $8931$.\\ 
	$ (16,4T,T^2 + 2T) $ & $ 4 $ & $ 2^6 $ & $4323$.\\ 
	$ (32,2T + 28,T^2 + 28) $ & $ 4 $ & $ 2^6 $ & $5183$.\\ 
	$ (16,8T,T^2 + 6T + 4) $ & $ 4 $ & $ 2^7 $ & $4827$, $6267$, $9523$.\\ 
	$ (16,4T,2T^2,T^3) $ & $ 4 $ & $ 2^7 $ & $6555$, $6699$, $8715$.\\ 
	$ (32,4T + 24,T^2 + 2T + 24) $ & $ 4 $ & $ 2^7 $ & $2811$, $7923$.\\ 
	$ (16,8T,T^2 + 2T + 12) $ & $ 4 $ & $ 2^7 $ & $1387$.\\ 
	$ (16,8T,T^2 + 14) $ & $ 4 $ & $ 2^7 $ & $8907$.\\ 
	$ (8,4T,2T^3 + 4,T^4 + 2T) $ & $ 4 $ & $ 2^8 $ & $1731$.\\ 
	$ (16,4T,2T^2,T^4) $ & $ 4 $ & $ 2^8 $ & $2415$.\\ 
	$ (16,8T,2T^2 + 4,T^3 + 2T) $ & $ 4 $ & $ 2^8 $ & $3059$.\\ 
	$ (32,2T + 12,T^4 + 16) $ & $ 4 $ & $ 2^8 $ & $3247$.\\ 
	$ (8,4T,2T^3 + 4,T^4 + 2T + 4) $ & $ 4 $ & $ 2^8 $ & $3855$.\\ 
	$ (16,8T,2T^2 + 8,T^3 + 4T) $ & $ 4 $ & $ 2^8 $ & $4899$.\\ 
	$ (16,T^2 + 6) $ & $ 4 $ & $ 2^8 $ & $5339$.\\ 
	$ (16,8T,2T^2 + 4T + 8,T^3 + 8) $ & $ 4 $ & $ 2^8 $ & $6083$.\\ 
	$ (16,T^2 + 12T + 14) $ & $ 4 $ & $ 2^8 $ & $9563$.\\ 
	$ (16,8T,2T^2 + 4T,T^3 + 4T) $ & $ 4 $ & $ 2^8 $ & $9835$.\\ 
	$ (16,2T^2 + 12,T^3 + 6T) $ & $ 4 $ & $ 2^9 $ & $3171$.\\ 
	$ (8,4T^2,2T^3 + 4T + 4,T^4 + 2T^2 + 2T + 4) $ & $ 4 $ & $ 2^9 $ & $7871$.\\ 
	$ (32,4T + 24,2T^2 + 24,T^4 + 16) $ & $ 4 $ & $ 2^9 $ & $8295$.\\ 
	$ (8,4T,2T^3 + 4,T^5 + 2T^2) $ & $ 4 $ & $ 2^9 $ & $8687$.\\ 
	$ (16,8T,2T^2 + 4T + 8,T^4) $ & $ 4 $ & $ 2^9 $ & $8979$.\\ 
	$ (16,2T^2 + 12,T^3 + 14T) $ & $ 4 $ & $ 2^9 $ & $9087$.\\ 
	$ (8,2T,T^8) $ & $ 4 $ & $ 2^{10} $ & $3007$.\\ 
	$ (8,2T^3 + 4T + 4,T^4 + 6T^2 + 2T) $ & $ 4 $ & $ 2^{10} $ & $5911$.\\ 
	$ (4,2T,T^9) $ & $ 4 $ & $ 2^{10} $ & $5983$.\\ 
	$ (16,8T,4T^2,2T^3,T^4) $ & $ 4 $ & $ 2^{10} $ & $8463$.\\ 
	$ (8,2T^3 + 4T^2 + 4T + 4,T^5 + 6T^2 + 4) $ & $ 4 $ & $ 2^{11} $ & $4879$.\\ 
	$ (8,4T,2T^3 + 4,T^{10} + 2T^2) $ & $ 4 $ & $ 2^{14} $ & $7967$.\\ 
	$ (32,2T + 24,T^2 + 16) $ & $ 5 $ & $ 2^6 $ & $4611$, $8347$.\\ 
	$ (32,2T + 8,T^2 + 16) $ & $ 5 $ & $ 2^6 $ & $6915$, $8723$.\\ 
	$ (32,2T + 16,T^2) $ & $ 5 $ & $ 2^6 $ & $7955$.\\ 
	$ (64,2T + 60,T^2 + 60) $ & $ 5 $ & $ 2^7 $ & $7051$.\\ 
	$ (32,4T + 16,2T^2,T^3) $ & $ 5 $ & $ 2^8 $ & $3243$, $7215$.\\ 
	$ (32,8T,T^2 + 2T + 8) $ & $ 5 $ & $ 2^8 $ & $9291$.\\ 
	$ (32,16T,T^2 + 6T + 4) $ & $ 5 $ & $ 2^9 $ & $627$.\\ 
	$ (32,16T,2T^2 + 8T + 16,T^3 + 8T) $ & $ 5 $ & $ 2^{10} $ & $4715$.\\ 
	$ (64,2T + 48,T^2) $ & $ 6 $ & $ 2^7 $ & $6095$.\\ 
	$ (64,2T + 8,T^2 + 48) $ & $ 6 $ & $ 2^7 $ & $7163$.\\ 
	$ (64,4T + 16,T^2 + 2T + 24) $ & $ 6 $ & $ 2^8 $ & $2451$.\\ 
	$ (64,4T + 48,T^2 + 2T + 40) $ & $ 6 $ & $ 2^8 $ & $6099$.\\ 
	$ (128,2T + 28,T^2 + 60) $ & $ 6 $ & $ 2^8 $ & $8299$.\\ 
	$ (64,8T + 32,T^2 + 2T + 56) $ & $ 6 $ & $ 2^9 $ & $2739$.\\ 
	$ (128,4T + 120,T^2 + 2T + 120) $ & $ 6 $ & $ 2^9 $ & $7491$.\\ 
	$ (64,4T + 16,2T^2 + 32,T^3) $ & $ 6 $ & $ 2^9 $ & $9503$.\\ 
	$ (64,16T,T^2 + 2T + 56) $ & $ 6 $ & $ 2^{10} $ & $8283$.\\ 
	$ (64,32T,T^2 + 24T + 54) $ & $ 6 $ & $ 2^{11} $ & $1851$.\\ 
	$ (128,2T + 104,T^2 + 112) $ & $ 7 $ & $ 2^8 $ & $1067$.\\ 
	$ (128,2T + 72,T^2 + 112) $ & $ 7 $ & $ 2^8 $ & $2595$.\\ 
	$ (128,2T + 112,T^2 + 64) $ & $ 7 $ & $ 2^8 $ & $5835$.\\ 
	$ (128,4T + 16,2T^2 + 96,T^3 + 64) $ & $ 7 $ & $ 2^{10} $ & $3003$.\\ 
	$ (256,128T,2T^2 + 44T + 68,T^3 + 62T + 20) $ & $ 7 $ & $ 2^{16} $ & $7535$.\\ 
	$ (512,2T + 28,T^2 + 316) $ & $ 8 $ & $ 2^{10} $ & $1023$.\\ 
	$ (512,2T + 304,T^2 + 448) $ & $ 9 $ & $ 2^{10} $ & $8679$.\\ 
	$ (2048,4T + 1520,T^2 + 2T + 744) $ & $ 11 $ & $ 2^{13} $ & $3363$.\\

\end{longtable}

\endgroup

\subsection*{Results for $f\equiv5\bmod8$} Let $3\le f<10000$ be a squarefree integer so  that $f$ is congruent to $5$ modulo $8$. Assume that the class number of $F_0=\Q(\sqrt{f})$ is even. As presented in Section \ref{algo}, we compute an ideal $J\subset\Z_2[X]$ such that for any positive $n$ we have a surjective morphism
\[ \frac{\Z_2[X]}{J}\longrightarrow A_n^{\perp}.\]
Denote by $N$ the index $[J:\Z_2[X]]$.
Let $n_0$ be the least integer such that the polynomial $\sum_{i=0}^{2^{n_0}-1}X^i$ lies in $J$. By Lemma \ref{relation-An}, the natural morphisms $A_{n_0}\longrightarrow A_n$ are isomorphisms for all $n\ge n_0$.  

We collect the ideals $J$ in the following table. It is convenient to write the results with the parameter $T=X-1$.

\begingroup
\fontsize{8}{10}\selectfont
\begin{longtable}{|p{0.33\textwidth}|p{0.02\textwidth}|p{0.02\textwidth}|p{0.52\textwidth}|}
	
	\hline
	$J$ & $n_0$ & $N$ & $f$ \\ 
	\hline
	\hline 
	\endhead
	
	\hline
	\endfoot
	
	$ (2,T^2) $ & $ 2 $ & $ 2^2 $ & $85$, $365$, $485$, $493$, $629$, $949$, $965$, $1037$, $1157$, $1165$, $1261$, $1853$, $2117$, $2165$, $2509$, $2581$, $2813$, $2941$, $3077$, $3085$, $3133$, $3293$, $3349$, $3365$, $3589$, $3653$, $3869$, $4573$, $4685$, $4709$, $4885$, $5165$, $5213$, $5389$, $5429$, $5485$, $5597$, $5837$, $5965$, $6485$, $6749$, $7141$, $7165$, $7373$, $8021$, $8485$, $8765$, $8917$, $8989$, $9197$, $9365$, $9565$, $9997$.\\ 
	$ (4,2T,T^2) $ & $ 2 $ & $ 2^3 $ & $165$, $205$, $221$, $285$, $429$, $445$, $741$, $885$, $901$, $1405$, $1517$, $1717$, $1749$, $1965$, $2013$, $2045$, $2085$, $2109$, $2245$, $2301$, $2365$, $2613$, $2685$, $2845$, $3005$, $3029$, $3165$, $3237$, $3597$, $3685$, $3765$, $3845$, $4069$, $4085$, $4453$, $4565$, $4645$, $4717$, $4773$, $4965$, $5045$, $5069$, $5133$, $5245$, $5629$, $5645$, $5685$, $5757$, $5829$, $5885$, $5917$, $5933$, $6149$, $6213$, $6285$, $6341$, $6357$, $6365$, $6445$, $6549$, $6605$, $6773$, $6837$, $6989$, $7045$, $7157$, $7221$, $7365$, $7405$, $7445$, $7485$, $7733$, $7813$, $7837$, $7869$, $7885$, $8333$, $8357$, $8493$, $8565$, $8749$, $8853$, $8877$, $8965$, $9005$, $9077$, $9141$, $9285$, $9309$, $9381$, $9469$, $9701$, $9885$.\\ 
	$ (4,2T,T^2 + 2) $ & $ 2 $ & $ 2^3 $ & $533$, $685$, $2285$, $4141$, $5317$, $6109$, $6437$, $7261$, $8285$, $9965$.\\ 
	$ (2,T^3) $ & $ 2 $ & $ 2^3 $ & $565$, $1285$, $1765$, $3277$, $4181$, $4589$, $5765$, $8885$, $9509$.\\ 
	$ (4,2T,T^3) $ & $ 2 $ & $ 2^4 $ & $861$, $1533$, $1645$, $1869$, $2093$, $2485$, $2821$, $3157$, $4277$, $4445$, $5285$, $5405$, $5453$, $5621$, $5781$, $5957$, $6141$, $6293$, $6573$, $6789$, $6853$, $6965$, $7805$, $8029$, $8165$, $8277$, $8533$, $8589$, $8733$, $9269$, $9373$, $9453$, $9541$, $9709$, $9717$.\\ 
	$ (8,2T + 4,T^2 + 4) $ & $ 2 $ & $ 2^4 $ & $645$, $957$, $1005$, $1205$, $3205$, $3333$, $3405$, $3741$, $4245$, $4605$, $5141$, $5709$, $6061$, $6501$, $6613$, $6645$, $7005$, $7845$, $8205$, $8621$, $9861$.\\ 
	$ (4,T^2 + 2T + 2) $ & $ 2 $ & $ 2^4 $ & $1565$, $2173$, $4469$, $4765$, $7421$, $9389$.\\ 
	$ (4,2T^2,T^3 + 2T) $ & $ 2 $ & $ 2^5 $ & $1085$, $2765$, $2877$, $3565$, $3813$, $4893$, $5901$, $5989$, $6685$, $7189$, $7453$, $9085$.\\ 
	$ (8,4T,T^2 + 2T + 6) $ & $ 2 $ & $ 2^5 $ & $1189$, $1781$, $4285$, $5941$, $9893$.\\ 
	$ (8,4T,T^2 + 2T + 2) $ & $ 2 $ & $ 2^5 $ & $3973$.\\ 
	$ (16,2T + 4,T^2 + 12) $ & $ 2 $ & $ 2^5 $ & $5205$.\\ 
	$ (8,4T,2T^2 + 4,T^3 + 2T + 4) $ & $ 2 $ & $ 2^6 $ & $2829$.\\ 
	$ (8,2T^2 + 4T + 4,T^3 + 6T + 4) $ & $ 2 $ & $ 2^7 $ & $9485$.\\ 
	$ (8,2T,T^2) $ & $ 3 $ & $ 2^4 $ & $1045$, $1677$, $2533$, $3245$, $5181$, $5973$, $6757$, $6981$, $7557$, $7645$, $7957$, $8149$, $8653$, $8949$, $9213$, $9789$, $9845$.\\ 
	$ (2,T^4) $ & $ 3 $ & $ 2^4 $ & $2885$, $2965$, $6085$, $7501$, $7709$.\\ 
	$ (4,2T,T^4) $ & $ 3 $ & $ 2^5 $ & $2037$, $2261$, $4381$, $4405$, $5061$, $5117$, $5253$, $6005$, $6693$, $7021$, $7429$, $7469$, $7701$, $7973$, $8421$, $9093$, $9429$, $9877$.\\ 
	$ (8,4T,T^2 + 4) $ & $ 3 $ & $ 2^5 $ & $2501$, $4045$, $5365$, $7085$, $7093$, $8077$, $9685$, $9805$.\\ 
	$ (8,4T,T^2) $ & $ 3 $ & $ 2^5 $ & $1885$, $3805$, $3965$, $7685$, $8845$.\\ 
	$ (16,2T + 12,T^2 + 12) $ & $ 3 $ & $ 2^5 $ & $2445$, $2717$, $3021$, $8805$, $9581$.\\ 
	$ (4,2T^2,T^3) $ & $ 3 $ & $ 2^5 $ & $3341$, $6245$, $6893$.\\ 
	$ (8,4T,2T^2,T^3) $ & $ 3 $ & $ 2^6 $ & $1365$, $3045$, $3885$, $4389$, $4485$, $5005$, $6405$, $7917$, $8645$, $9165$, $9597$, $9933$.\\ 
	$ (4,2T^2,T^4 + 2T) $ & $ 3 $ & $ 2^6 $ & $357$, $1173$, $1309$, $3621$, $4301$, $8517$.\\ 
	$ (4,2T,T^5) $ & $ 3 $ & $ 2^6 $ & $1581$, $2373$, $4029$, $5797$, $7797$.\\ 
	$ (4,2T^2,T^4) $ & $ 3 $ & $ 2^6 $ & $3485$, $6205$, $9061$.\\ 
	$ (8,4T,2T^2 + 4,T^3 + 2T) $ & $ 3 $ & $ 2^6 $ & $5845$, $9821$.\\ 
	$ (4,2T^2,T^4 + 2) $ & $ 3 $ & $ 2^6 $ & $7765$, $9773$.\\ 
	$ (16,4T + 8,T^2 + 12) $ & $ 3 $ & $ 2^6 $ & $2405$.\\ 
	$ (32,2T + 20,T^2 + 28) $ & $ 3 $ & $ 2^6 $ & $3893$.\\ 
	$ (4,2T^2,T^5 + 2T) $ & $ 3 $ & $ 2^7 $ & $6477$, $7285$, $8789$, $9741$.\\ 
	$ (8,2T^2 + 4,T^3 + 6T) $ & $ 3 $ & $ 2^7 $ & $805$, $1469$, $3605$.\\ 
	$ (8,4T,2T^2,T^4) $ & $ 3 $ & $ 2^7 $ & $2805$, $4053$.\\ 
	$ (8,2T^2 + 4T + 4,T^3 + 2T + 4) $ & $ 3 $ & $ 2^7 $ & $6461$, $8365$.\\ 
	$ (16,4T + 8,2T^2 + 8,T^3 + 8) $ & $ 3 $ & $ 2^7 $ & $5037$.\\ 
	$ (16,4T + 8,2T^2,T^3) $ & $ 3 $ & $ 2^7 $ & $8005$.\\ 
	$ (4,2T^2,T^6) $ & $ 3 $ & $ 2^8 $ & $8245$, $9605$.\\ 
	$ (4,T^4 + 2T^3 + 2) $ & $ 3 $ & $ 2^8 $ & $1685$.\\ 
	$ (8,4T,2T^3,T^4 + 2T^2) $ & $ 3 $ & $ 2^8 $ & $4845$.\\ 
	$ (4,2T^3,T^5 + 2T^2 + 2T) $ & $ 3 $ & $ 2^8 $ & $5397$.\\ 
	$ (8,4T,2T^3,T^4) $ & $ 3 $ & $ 2^8 $ & $6765$.\\ 
	$ (16,4T + 8,2T^2 + 8,T^4) $ & $ 3 $ & $ 2^8 $ & $7565$.\\ 
	$ (8,4T,2T^2,T^5) $ & $ 3 $ & $ 2^8 $ & $8701$.\\ 
	$ (16,4T + 8,2T^2 + 8,T^5) $ & $ 3 $ & $ 2^9 $ & $7077$.\\ 
	$ (4,2T^4,T^5 + 2T) $ & $ 3 $ & $ 2^9 $ & $7413$.\\ 
	$ (16,2T + 8,T^2) $ & $ 4 $ & $ 2^5 $ & $2005$, $4981$, $5109$, $5605$, $7437$, $8605$.\\ 
	$ (32,2T + 12,T^2 + 28) $ & $ 4 $ & $ 2^6 $ & $2701$.\\ 
	$ (32,2T + 28,T^2 + 28) $ & $ 4 $ & $ 2^6 $ & $9645$.\\ 
	$ (64,2T + 52,T^2 + 28) $ & $ 4 $ & $ 2^7 $ & $1653$, $4917$, $8445$.\\ 
	$ (16,8T,T^2 + 4T + 8) $ & $ 4 $ & $ 2^7 $ & $2605$.\\ 
	$ (16,8T,T^2 + 4T + 12) $ & $ 4 $ & $ 2^7 $ & $3445$.\\ 
	$ (16,8T,2T^2,T^3) $ & $ 4 $ & $ 2^8 $ & $4669$.\\ 
	$ (32,4T + 24,2T^2 + 24,T^3 + 24) $ & $ 4 $ & $ 2^8 $ & $5565$.\\ 
	$ (32,8T + 16,T^2 + 4) $ & $ 4 $ & $ 2^8 $ & $7397$.\\ 
	$ (16,2T^2 + 12T + 12,T^3 + 10T + 12) $ & $ 4 $ & $ 2^9 $ & $9205$.\\ 
	$ (16,8T,2T^3 + 4T^2 + 8,T^4 + 4T + 8) $ & $ 4 $ & $ 2^{11} $ & $7293$.\\ 
	$ (8,2T^4 + 4T + 4,T^5 + 6T^2 + 6T) $ & $ 4 $ & $ 2^{13} $ & $2397$.\\ 
	$ (16,4T,2T^3,T^9 + 2T^2) $ & $ 4 $ & $ 2^{14} $ & $9021$.\\ 
	$ (32,2T + 8,T^2 + 16) $ & $ 5 $ & $ 2^6 $ & $1221$, $3477$.\\ 
	$ (32,2T + 24,T^2 + 16) $ & $ 5 $ & $ 2^6 $ & $7205$, $9797$.\\ 
	$ (32,2T + 16,T^2) $ & $ 5 $ & $ 2^6 $ & $8229$.\\ 
	$ (32,8T,T^2 + 8) $ & $ 5 $ & $ 2^8 $ & $6565$.\\ 
	$ (128,2T + 20,T^2 + 28) $ & $ 5 $ & $ 2^8 $ & $8437$.\\ 
	$ (64,4T + 56,2T^2 + 56,T^3 + 56) $ & $ 5 $ & $ 2^9 $ & $6045$.\\ 
	$ (16,8T,2T^4 + 4T^3 + 4T^2 + 4T + 12,T^5 + 6T^3 + 6T^2 + 2T + 4) $ & $ 5 $ & $ 2^{14} $ & $9445$.\\ 
	$ (64,2T + 56,T^2 + 48) $ & $ 6 $ & $ 2^7 $ & $2669$.\\ 
	$ (64,2T + 8,T^2 + 48) $ & $ 6 $ & $ 2^7 $ & $4173$.\\ 
	$ (128,2T + 92,T^2 + 60) $ & $ 6 $ & $ 2^8 $ & $1245$.\\ 
	$ (128,2T + 76,T^2 + 92) $ & $ 6 $ & $ 2^8 $ & $9669$.\\ 
	$ (128,2T + 120,T^2 + 112) $ & $ 7 $ & $ 2^8 $ & $5421$.\\ 
	$ (256,2T + 72,T^2 + 240) $ & $ 8 $ & $ 2^9 $ & $6805$.\\ 
	$ (1024,2T + 316,T^2 + 636) $ & $ 9 $ & $ 2^{11} $ & $1605$.\\ 
	$ (2048,2T + 1036,T^2 + 2012) $ & $ 10 $ & $ 2^{12} $ & $8045$.\\ 
	
\end{longtable}
\endgroup

\subsection*{Results for $f\equiv1\bmod8$} Let $3\le f<10000$ be a squarefree integer so  that $f$ is congruent to $1$ modulo 8. As presented in Section \ref{algo}, we compute an ideal $J\subset\Z_2[X]$ such that for any positive $n$ we have a surjective morphism
\[ \frac{\Z_2[X]}{J}\longrightarrow \left(\frac{A_n}{A_0}\right)^{\perp}.\]
Denote by $N$ the index $[\Z_2[X]:J]$.
Let $n_0$ be the least integer such that the polynomial $\sum_{i=0}^{2^{n_0}-1}X^i$ lies in $J$. By Lemma \ref{relation-An}, the natural morphisms $A_{n_0}\longrightarrow A_n$ are isomorphisms for all $n\ge n_0$.  

We collect the ideals $J$ in the following table. It is convenient to write the results with the parameter $T=X-1$. We omit the integers $f$ such that the ideal $J=\Z_2[T]$.

\begingroup
\fontsize{8}{10}\selectfont
\begin{longtable}{|p{0.33\textwidth}|p{0.02\textwidth}|p{0.02\textwidth}|p{0.52\textwidth}|}
	
	\hline
	$J$ & $n_0$ & $N$ & $f$ \\ 
	\hline
	\hline 
	\endhead
	
	\hline
	\endfoot
	
	$ (2,T) $ & $ 1 $ & $ 2 $ & $105$, $145$, $185$, $217$, $273$, $329$, $345$, $385$, $409$, $481$, $521$, $553$, $569$, $665$, $705$, $713$, $785$, $809$, $857$, $865$, $897$, $953$, $1001$, $1129$, $1145$, $1209$, $1265$, $1281$, $1337$, $1385$, $1417$, $1465$, $1505$, $1545$, $1657$, $1729$, $1745$, $1769$, $1817$, $1833$, $1865$, $1905$, $1937$, $2065$, $2121$, $2137$, $2153$, $2185$, $2249$, $2289$, $2345$, $2377$, $2505$, $2521$, $2553$, $2561$, $2585$, $2617$, $2633$, $2681$, $2697$, $2705$, $2729$, $2769$, $2785$, $2905$, $3081$, $3161$, $3209$, $3241$, $3297$, $3337$, $3345$, $3353$, $3441$, $3601$, $3633$, $3657$, $3665$, $3745$, $3769$, $3785$, $3801$, $3809$, $3905$, $3913$, $3945$, $4081$, $4089$, $4209$, $4345$, $4385$, $4393$, $4433$, $4465$, $4553$, $4585$, $4697$, $4809$, $4865$, $4929$, $4945$, $4953$, $5033$, $5129$, $5177$, $5305$, $5353$, $5369$, $5497$, $5505$, $5545$, $5585$, $5609$, $5649$, $5665$, $5673$, $5705$, $5745$, $5777$, $5801$, $5817$, $5889$, $5897$, $5993$, $6041$, $6097$, $6153$, $6161$, $6169$, $6217$, $6265$, $6329$, $6385$, $6513$, $6649$, $6721$, $6785$, $6873$, $6937$, $6945$, $6985$, $7033$, $7097$, $7177$, $7193$, $7241$, $7273$, $7305$, $7329$, $7385$, $7449$, $7505$, $7521$, $7545$, $7553$, $7657$, $7705$, $7761$, $7777$, $7833$, $7849$, $7945$, $7969$, $8009$, $8057$, $8065$, $8113$, $8137$, $8185$, $8305$, $8321$, $8329$, $8337$, $8377$, $8441$, $8465$, $8497$, $8617$, $8665$, $8697$, $8769$, $8785$, $8841$, $8953$, $9017$, $9089$, $9105$, $9145$, $9161$, $9169$, $9185$, $9193$, $9289$, $9321$, $9353$, $9361$, $9433$, $9545$, $9577$, $9593$, $9665$, $9681$, $9705$, $9737$, $9785$, $9817$, $9865$, $9889$, $9905$, $9913$, $9985$.\\ 
	$ (4,T + 2) $ & $ 1 $ & $ 2^2 $ & $113$, $497$, $1601$, $1777$, $1841$, $2369$, $3073$, $3121$, $3137$, $3473$, $4177$, $4193$, $4721$, $4817$, $5761$, $6353$, $6481$, $7153$, $7313$, $7489$, $8209$, $9569$, $9649$.\\ 
	$ (8,T + 2) $ & $ 1 $ & $ 2^3 $ & $3089$, $3313$, $3521$, $4513$, $9793$.\\ 
	$ (4,T) $ & $ 2 $ & $ 2^2 $ & $65$, $137$, $265$, $777$, $985$, $1065$, $1073$, $1313$, $1321$, $1673$, $1705$, $1897$, $2001$, $2233$, $2257$, $2265$, $2297$, $2713$, $2777$, $2921$, $2929$, $2945$, $2985$, $3001$, $3193$, $3265$, $3585$, $3857$, $3865$, $3985$, $4033$, $4065$, $4121$, $4137$, $4649$, $4665$, $4889$, $5065$, $5161$, $5217$, $5257$, $5273$, $5385$, $5417$, $5449$, $5465$, $5473$, $5681$, $6065$, $6145$, $6177$, $6233$, $6377$, $6465$, $6473$, $6585$, $6641$, $6665$, $6745$, $6865$, $6969$, $7001$, $7049$, $7185$, $7337$, $7369$, $7417$, $7433$, $7745$, $7881$, $7897$, $8153$, $8345$, $8593$, $8601$, $8681$, $8961$, $8985$, $9393$, $9689$, $9833$, $9961$.\\ 
	$ (2,T^2) $ & $ 2 $ & $ 2^2 $ & $697$, $969$, $2193$, $3009$, $3417$, $3553$, $3649$, $3977$, $4233$, $4521$, $5289$, $5321$, $5529$, $6369$, $7257$, $7809$, $7913$, $7953$, $8041$, $8585$, $8857$, $9129$, $9553$.\\ 
	$ (2,T^3) $ & $ 2 $ & $ 2^3 $ & $337$, $1457$, $1553$, $2129$, $2449$, $3729$, $5233$, $5969$, $7121$, $7409$, $8249$, $8401$, $8481$, $8977$.\\ 
	$ (4,2T,T^2) $ & $ 2 $ & $ 2^3 $ & $1241$, $1513$, $2993$, $4161$, $4505$, $4785$, $5785$, $6105$, $6953$, $8177$, $8265$, $9673$.\\ 
	$ (4,2T,T^2 + 2) $ & $ 2 $ & $ 2^3 $ & $593$, $881$, $2689$, $4049$, $5945$, $6689$, $9137$, $9417$, $9697$.\\ 
	$ (8,T + 6) $ & $ 2 $ & $ 2^3 $ & $2593$, $4417$, $5297$, $6049$, $8561$.\\ 
	$ (4,2T,T^3) $ & $ 2 $ & $ 2^4 $ & $6601$, $7345$, $7665$, $8905$, $9177$, $9345$.\\ 
	$ (4,T^2 + 2) $ & $ 2 $ & $ 2^4 $ & $577$, $2409$, $2937$.\\ 
	$ (4,T^2 + 2T + 2) $ & $ 2 $ & $ 2^4 $ & $4001$, $4657$.\\ 
	$ (8,2T + 4,T^2) $ & $ 2 $ & $ 2^4 $ & $8569$.\\ 
	$ (8,2T + 4,T^2 + 4) $ & $ 2 $ & $ 2^4 $ & $9881$.\\ 
	$ (8,2T,T^3 + 4) $ & $ 2 $ & $ 2^5 $ & $9377$.\\ 
	$ (8,4T,2T^2 + 4,T^3 + 6T + 4) $ & $ 2 $ & $ 2^6 $ & $2665$.\\ 
	$ (8,2T^2 + 4T + 4,T^3 + 6T + 4) $ & $ 2 $ & $ 2^7 $ & $5313$.\\ 
	$ (8,T) $ & $ 3 $ & $ 2^3 $ & $41$, $313$, $457$, $761$, $1081$, $1561$, $1985$, $1993$, $2305$, $2545$, $2569$, $2849$, $2953$, $2977$, $3017$, $3129$, $3385$, $3929$, $4017$, $4249$, $4681$, $4729$, $4793$, $4841$, $4921$, $5641$, $5713$, $6073$, $6185$, $6401$, $6617$, $6657$, $7465$, $7473$, $7561$, $7721$, $7801$, $8169$, $8393$, $8521$, $8545$, $8729$, $8809$, $9049$, $9257$, $9305$, $9385$, $9465$, $9529$, $9929$.\\ 
	$ (8,T + 4) $ & $ 3 $ & $ 2^3 $ & $257$, $1153$, $2513$, $4289$, $5201$, $6209$, $6449$, $6529$, $8081$.\\ 
	$ (8,2T,T^2) $ & $ 3 $ & $ 2^4 $ & $1105$, $3201$, $6409$.\\ 
	$ (4,T^2) $ & $ 3 $ & $ 2^4 $ & $561$, $6681$.\\ 
	$ (4,T^2 + 2T) $ & $ 3 $ & $ 2^4 $ & $1353$, $8313$.\\ 
	$ (2,T^4) $ & $ 3 $ & $ 2^4 $ & $1921$, $5729$.\\ 
	$ (16,T + 14) $ & $ 3 $ & $ 2^4 $ & $6433$, $7217$.\\ 
	$ (16,T + 6) $ & $ 3 $ & $ 2^4 $ & $5393$.\\ 
	$ (4,2T,T^3 + 2) $ & $ 3 $ & $ 2^4 $ & $8161$.\\ 
	$ (8,2T,T^2 + 4) $ & $ 3 $ & $ 2^4 $ & $8273$.\\ 
	$ (4,2T,T^4) $ & $ 3 $ & $ 2^5 $ & $1785$, $5865$.\\ 
	$ (8,4T,T^2 + 6) $ & $ 3 $ & $ 2^5 $ & $4273$.\\ 
	$ (4,2T^2,T^3) $ & $ 3 $ & $ 2^5 $ & $5617$.\\ 
	$ (8,4T,T^2 + 2T) $ & $ 3 $ & $ 2^5 $ & $8633$.\\ 
	$ (16,2T + 12,T^2 + 12) $ & $ 3 $ & $ 2^5 $ & $8745$.\\ 
	$ (4,2T,T^5) $ & $ 3 $ & $ 2^6 $ & $3689$, $7905$.\\ 
	$ (8,2T + 4,T^4) $ & $ 3 $ & $ 2^6 $ & $3281$.\\ 
	$ (8,2T,T^4) $ & $ 3 $ & $ 2^6 $ & $6545$.\\ 
	$ (16,4T + 8,T^2 + 8) $ & $ 3 $ & $ 2^6 $ & $7089$.\\ 
	$ (16,4T + 8,T^2 + 2T + 4) $ & $ 3 $ & $ 2^6 $ & $7689$.\\ 
	$ (4,2T^3,T^4) $ & $ 3 $ & $ 2^7 $ & $1649$, $2737$, $8449$.\\ 
	$ (2,T^7) $ & $ 3 $ & $ 2^7 $ & $3937$, $5921$.\\ 
	$ (4,2T^2,T^5 + 2T) $ & $ 3 $ & $ 2^7 $ & $4369$.\\ 
	$ (8,4T,2T^2,T^4 + 2T) $ & $ 3 $ & $ 2^7 $ & $4641$.\\ 
	$ (8,4T,2T^2,T^4 + 2T + 4) $ & $ 3 $ & $ 2^7 $ & $7361$.\\ 
	$ (4,2T^3,T^4 + 2T^2) $ & $ 3 $ & $ 2^7 $ & $7633$.\\ 
	$ (8,2T,T^5) $ & $ 3 $ & $ 2^7 $ & $9401$.\\ 
	$ (8,4T,2T^3,T^5 + 2T^2 + 2T) $ & $ 3 $ & $ 2^9 $ & $5593$.\\ 
	$ (16,T) $ & $ 4 $ & $ 2^4 $ & $465$, $609$, $889$, $1585$, $2865$, $3065$, $3233$, $3593$, $4073$, $4537$, $4969$, $4985$, $5849$, $7985$, $8089$, $8945$, $9217$, $9641$, $9953$.\\ 
	$ (16,T + 8) $ & $ 4 $ & $ 2^4 $ & $1169$, $1393$, $1633$, $4529$, $9121$.\\ 
	$ (16,T + 12) $ & $ 4 $ & $ 2^4 $ & $4577$, $6881$, $7393$.\\ 
	$ (16,T + 4) $ & $ 4 $ & $ 2^4 $ & $5089$.\\ 
	$ (16,2T + 8,T^2) $ & $ 4 $ & $ 2^5 $ & $3145$, $5185$.\\ 
	$ (16,2T,T^2) $ & $ 4 $ & $ 2^5 $ & $3961$, $4777$.\\ 
	$ (16,2T,T^2 + 8) $ & $ 4 $ & $ 2^5 $ & $2113$.\\ 
	$ (32,T + 14) $ & $ 4 $ & $ 2^5 $ & $7441$.\\ 
	$ (8,T^2 + 6T + 4) $ & $ 4 $ & $ 2^6 $ & $5073$, $9473$.\\ 
	$ (16,2T + 8,T^3 + 8) $ & $ 4 $ & $ 2^6 $ & $3761$.\\ 
	$ (16,2T,T^3) $ & $ 4 $ & $ 2^6 $ & $4305$.\\ 
	$ (32,2T + 28,T^2 + 28) $ & $ 4 $ & $ 2^6 $ & $4745$.\\ 
	$ (8,T^2 + 4T) $ & $ 4 $ & $ 2^6 $ & $5457$.\\ 
	$ (32,8T + 16,T^2 + 24) $ & $ 4 $ & $ 2^8 $ & $1217$.\\ 
	$ (16,4T + 8,T^3 + 2T^2 + 2T) $ & $ 4 $ & $ 2^8 $ & $6441$.\\ 
	$ (32,8T + 16,T^2 + 4T + 24) $ & $ 4 $ & $ 2^8 $ & $9281$.\\ 
	$ (16,8T,4T^2,T^3 + 2T^2 + 2T) $ & $ 4 $ & $ 2^9 $ & $4633$.\\ 
	$ (16,2T^2 + 4,T^3 + 2T) $ & $ 4 $ & $ 2^9 $ & $7161$.\\ 
	$ (4,2T^3,T^7 + 2T^2 + 2T + 2) $ & $ 4 $ & $ 2^{10} $ & $6913$.\\ 
	$ (8,4T,2T^3,T^6 + 2T) $ & $ 4 $ & $ 2^{10} $ & $9809$.\\ 
	$ (8,4T^2,T^5 + 2T^3 + 6T + 6) $ & $ 4 $ & $ 2^{12} $ & $4801$.\\ 
	$ (32,T) $ & $ 5 $ & $ 2^5 $ & $1113$, $3289$, $3433$, $3545$, $5657$, $7321$, $7673$, $8297$, $9497$.\\ 
	$ (32,T + 8) $ & $ 5 $ & $ 2^5 $ & $161$, $2273$.\\ 
	$ (32,T + 16) $ & $ 5 $ & $ 2^5 $ & $2657$, $2833$.\\ 
	$ (32,T + 24) $ & $ 5 $ & $ 2^5 $ & $353$.\\ 
	$ (32,T + 20) $ & $ 5 $ & $ 2^5 $ & $1249$.\\ 
	$ (32,T + 4) $ & $ 5 $ & $ 2^5 $ & $8257$.\\ 
	$ (32,2T + 16,T^2) $ & $ 5 $ & $ 2^6 $ & $3705$.\\ 
	$ (32,2T,T^2 + 16) $ & $ 5 $ & $ 2^6 $ & $9273$.\\ 
	$ (32,4T,T^2 + 2T) $ & $ 5 $ & $ 2^7 $ & $7769$.\\ 
	$ (64,2T + 44,T^2 + 28) $ & $ 5 $ & $ 2^7 $ & $9265$.\\ 
	$ (32,4T,T^3 + 2T^2 + 2T + 24) $ & $ 5 $ & $ 2^9 $ & $7841$.\\ 
	$ (64,8T + 48,T^2 + 6T) $ & $ 5 $ & $ 2^9 $ & $8385$.\\ 
	$ (16,8T,T^3 + 4T^2 + 4T + 2) $ & $ 5 $ & $ 2^{10} $ & $1889$.\\ 
	$ (16,4T^2,T^3 + 12T + 12) $ & $ 5 $ & $ 2^{10} $ & $9601$.\\ 
	$ (64,T + 32) $ & $ 6 $ & $ 2^6 $ & $721$, $1057$, $5569$, $8369$.\\ 
	$ (64,T) $ & $ 6 $ & $ 2^6 $ & $2201$, $4441$, $7265$, $7609$.\\ 
	$ (64,T + 24) $ & $ 6 $ & $ 2^6 $ & $3409$.\\ 
	$ (128,T + 70) $ & $ 6 $ & $ 2^7 $ & $3361$.\\ 
	$ (64,2T + 16,T^2) $ & $ 6 $ & $ 2^7 $ & $6497$.\\ 
	$ (64,2T,T^2 + 32) $ & $ 6 $ & $ 2^7 $ & $7793$.\\ 
	$ (128,T + 30) $ & $ 6 $ & $ 2^7 $ & $9233$.\\ 
	$ (256,T + 202) $ & $ 6 $ & $ 2^8 $ & $2177$.\\ 
	$ (64,4T + 32,T^2 + 48) $ & $ 6 $ & $ 2^8 $ & $8609$.\\ 
	$ (64,2T + 16,T^3 + 32) $ & $ 6 $ & $ 2^8 $ & $9521$.\\ 
	$ (64,4T,2T^2,T^4 + 32) $ & $ 6 $ & $ 2^{10} $ & $6817$.\\ 
	$ (64,16T,2T^2 + 48,T^3) $ & $ 6 $ & $ 2^{11} $ & $8897$.\\ 
	$ (64,8T,4T^2,T^5 + 2T^3 + 2T) $ & $ 6 $ & $ 2^{15} $ & $6001$.\\ 
	$ (128,T + 16) $ & $ 7 $ & $ 2^7 $ & $3841$.\\ 
	$ (128,T) $ & $ 7 $ & $ 2^7 $ & $4265$.\\ 
	$ (128,T + 64) $ & $ 7 $ & $ 2^7 $ & $5873$.\\ 
	$ (128,2T + 96,T^2) $ & $ 7 $ & $ 2^8 $ & $2145$, $2329$.\\ 
	$ (128,2T + 32,T^2) $ & $ 7 $ & $ 2^8 $ & $2465$, $6305$.\\ 
	$ (128,2T,T^2) $ & $ 7 $ & $ 2^8 $ & $7081$.\\ 
	$ (128,2T + 32,T^2 + 64) $ & $ 7 $ & $ 2^8 $ & $8241$.\\ 
	$ (256,T + 64) $ & $ 8 $ & $ 2^8 $ & $3217$.\\ 
	$ (256,T + 160) $ & $ 8 $ & $ 2^8 $ & $6769$.\\ 
	$ (256,16T,8T^3,T^4 + 6T^3 + 14T^2 + 6T + 224) $ & $ 8 $ & $ 2^{19} $ & $4097$.\\ 
	$ (512,T) $ & $ 9 $ & $ 2^9 $ & $1185$.\\ 
	$ (512,T + 256) $ & $ 9 $ & $ 2^9 $ & $4993$.\\ 
	$ (512,2T + 176,T^2 + 192) $ & $ 9 $ & $ 2^{10} $ & $2337$.\\ 
	$ (512,4T,2T^2,T^3) $ & $ 9 $ & $ 2^{12} $ & $7585$.\\ 
	$ (2048,T) $ & $ 11 $ & $ 2^{11} $ & $3881$.\\ 
	$ (2048,2T,T^2 + 1024) $ & $ 11 $ & $ 2^{12} $ & $1201$.\\ 
	$ (2048,8T + 256,T^2 + 4T + 384) $ & $ 11 $ & $ 2^{14} $ & $4481$.\\ 
	$ (2048,8T + 256,2T^2 + 1536,T^3 + 4T + 128) $ & $ 11 $ & $ 2^{15} $ & $3713$.\\
	
\end{longtable}
\endgroup

Finally, we compare our results with the ones known in the literature, we denote by $p,q,r$ prime numbers. In \cite{ozakitaya} it has been shown that Greenberg's conjecture is true in the following cases:
\begin{enumerate}[label = (\roman*)]
	\item\label{ot1} If $f=p$ with $p\equiv1\bmod8$ and $\big(\frac{2}{p}\big)_4\left(\frac{p}{2}\right)_4=-1$ where $\left(\frac{*}{*}\right)_4$ denotes the bi-quadratic residue symbol, then $\A_{F_n}$ is trivial for all $n$.
	\item\label{ot2} If $f=pq$ with $p\equiv3$, $q\equiv3\bmod 8$, then $\A_n$ is trivial for all $n$.
	\item\label{ot3} If $f=pq$ with $p\equiv3\bmod 4$ and $q\equiv5\bmod 8$, then $\A_0\cong\Z/2\Z$ and Greenberg's conjecture holds.
	\item\label{ot4} If $f=pq$ with $p\equiv5$ and $q\equiv5\bmod 8$, then Greenberg's holds.	
\end{enumerate}

With our algorithms, in cases \ref{ot1} and \ref{ot2} we find the ideals $\Z_2[T]$ for all $f<10000$ getting that $A_n$, and hence $\A_n$ stabilize at level $0$. In case \ref{ot3} we find the ideal $(2,T)$ for all $f<10000$, therefore we get that $\#\A_1\le2\#A_1=4$ and the stabilization occurs before level $n=1$. However, we discover with Pari-GP that $\#\A_1=2$. In the last case, we find ideals of the form $(2^d,T)$ for some positive integer $d$, getting that $A_n$ and $\A_n$ are stable from level $d$.  

In \cite{mizusawa} the case where $f=pqr$ such that $p\equiv5$, $q\equiv5$ and $r\equiv3\bmod8$ and $\left(\frac{pq}{r}\right)=-1$ is studied. For such $f$ it is shown that $\#\A_n=\#\A_0=4$ for any $n\ge0$. Under these assumptions the algorithm computes the ideal $(4,2T,T^2)$ implying the stabilization before level $2$ and the relation $\#\A_n\le 2\#A_2=16$.

In \cite{fukudakomatsuprev} it is shown that if we assume $f=pq$ with $p\equiv3$ and $q\equiv1\bmod8$ and $\left(\frac{p}{q}\right)=-1$ and $2^{\frac{q-1}{4}}\equiv-1\bmod q$, then $\A_n\cong\A_1$ for any $n\ge1$. The algorithm computes the ideals $(2,T^2)$ getting that the stabilization occurs before level $2$. 

Finally, we assume that $f=pq$ such that $p\equiv7$ and $q\equiv1\bmod8$ and $\left(\frac{p}{q}\right)=-1$. This case is studied in \cite{mizusawa}; if $q\equiv1\bmod16$ or if $q\equiv7\bmod16$ and $2^{\frac{q-1}{4}}\equiv-1\bmod q$, then $\A_n\cong\A_1$. We  find the ideals $(2,T^2)$ in the first case and $(2,T^3)$ in the second one. Therefore, we get that the stabilization occurs before level $2$ or $3$ respectively. 

Notice that in all these cases our algorithm is less precise than the above mentioned results largely because of the assumption we explained at the beginning of Section \ref{furtherrmk}. 

\vspace{1.5cm}
\appendix
\section{Sinnott's index formula}\label{appA}

In the appendix we recall the definition of the group of cyclotomic units following \cite[Section 4]{sinnott}. This group is crucial in order to define the group $C_{F_n}$ in Definition \ref{new-group-of-cyclounits}. We denote by $\ze{n}$ the primitive $n$-th root of unity $\mathrm{e}^{\frac{2\pi i}{n}}$ for any integer $n\ge1$.

\begin{definition}
	Let $K$ be an abelian number field. For any integer $n\ge 1$ we denote by $K_{(n)}$ the intersection $K\cap\Q(\ze{n})$. Let $D_K$ be the subgroup of $K^\times$ generated by $-1$ and the elements
	\[ \underset{\Q(\ze{n})|K_{(n)}}{\norm}(1-\ze{n}^a) \quad:\quad n\in\Z_{>2},\,a\in\left\lbrace 1,2,\dots,n-1\right\rbrace. \]
	The group of cyclotomic units is $\cyc_K=D_K\cap O_K^\times$. 
\end{definition}
\begin{remark}\label{cyclo-DK}
	Let $G=\gal(K|\Q)$. It is easy to show that
	\[ D_K=\Q^\times\cdot\left\langle \underset{\Q(\ze{n})|K_{(n)}}{\norm}(1-\ze{n})\right\rangle_{\Z[G]}  \]
	where $n$ runs over the conductors of all the subfields of $K$.
\end{remark}

We focus on $K=F_n$ the $n$-th level field in the cyclotomic $\Z_2$-extension of the field $F=\Q(\sqrt{f})$. Sinnott's index formula \cite[Theorem 4.1]{sinnott} states that
\begin{equation}\label{index-sinnott}\tag*{(S)}
	[O_{F_n}^\times:\cyc_{F_n}]=\begin{cases}
		2^{2^{n+1}-1}\#\cl_{F_n}c_{F_n}, & \t{if $f\equiv1\bmod4$ and $f$ prime,}\\
		2^{2^{n+1}-2}\#\cl_{F_n}c_{F_n}, & \t{otherwise.}
	\end{cases}
\end{equation} 
The constant $c_{F_n}$ is a positive rational number.
In this appendix we give the definition of $c_{F_n}$ and we prove that $c_{F_n}=1$.

We introduce the notation used in \cite{sinnott}. Let $\Gamma_n=\gal(F_n|\Q)$.
\begin{itemize}
	\item For any prime number $p$, let $T_p<\Gamma_n$ be the inertia subgroup at $p$. Let $r$ be a positive squarefree integer. Denote by $T_r$ the subgroup of $\Gamma_n$ generated by $T_p$ for any prime $p$ that divides $r$.
	\item For any subgroup $H<\Gamma_n$ we put
	\[ s(H)=\sum_{h\in H}h. \]
	\item Let $p$ be a prime number. We define $(p,F_n)^*\in\Q[\Gamma_n]$ as 
	\[ (p,F_n)^*=\frac{s(T_p)}{\#T_p}\lambda_p^{-1} \]
	where $\lambda_p$ is a Frobenius element above $p$. Namely, for a chosen prime $\mathfrak{p}$ that divides $p$ then $\lambda_p(x)\equiv x^p\bmod\mathfrak{p}$. Since $F_n$ is an abelian field, $\lambda_p$ does not depends on the choice of $\mathfrak{p}$. Finally, notice that $\lambda_p$ is defined up to an element of $T_p$, but $(p,F_n)^*$ is well defined.
	\item Let $R=\Z[\Gamma_n]$. Let $m$ be the product of distinct primes that divide the conductor of $F_n$. Define $U\subset\Q[\Gamma_n]$ as the $R$-module generated by
	\[ s(T_r)\prod_{p|\frac{m}{r}}(1-(p,F_n)^*) \]
	where $r$ runs over the positive divisors of $m$ and $p$ denotes a prime number.
	Then \cite[Proposition $2.3$]{sinnott} states that $U$ is a lattice in $\Q(\Gamma_n)$.
\end{itemize}

Let $V$ be a finite dimensional $\Q$-vector space. Let $L,M\subset V$ be two lattices. Choose a $\Q$-linear automorphism $\phi$ of $V$ such that $\phi(L)=M$. Then, we define
\[ (L:M)=|\det(\phi)|. \]
The index $(L:M)$ does not depend on the choice of $\phi$. 
By \cite[Lemmas $1.1$ and $1.2$]{sinnott} the index satisfies the following properties.
Let $L,M,N\subset V$ be lattices then
\begin{equation}\label{index1} \tag*{($A1$)}
	(L:N)=(L:M)(M:N).
\end{equation}  
Moreover, let $L,M\subset V$ be both lattices and $R$-modules. For any $\alpha\in\Q(\Gamma_n)$ let $L[\alpha]$ be the kernel of the map $L\longrightarrow\alpha L$ induced by the multiplication by $\alpha$. Define $M[\alpha]$ similarly. Then
\begin{equation}\label{index2} \tag*{($A2$)}
	(L:M)=(L[\alpha]:M[\alpha])(\alpha L:\alpha M).
\end{equation}	
Finally, according to \cite[Theorem $4.1$]{sinnott} we have $c_{F_n}=(R:U)$.

\begin{proposition}
	Let $f\ge3$ be an odd squarefree integer. Let $F_n$ be the $n$-th level in the cyclotomic $\Z_2$-extension of $F=\Q(\sqrt{f})$. Then $c_{F_n}=1$. 
\end{proposition}

\begin{proof}
	First of all, we assume $n=0$. Since the group $\gal(F_0|\Q)$ is cyclic, the proposition follows from \cite[Theorem $5.4$]{sinnott}. 
	
	Fix a level $n\ge1$. Define $U_f\subset\Q(\Gamma_n)$ be the $R$-module generated by
	\[ s(T_r)\prod_{p|\frac{f}{r}}(1-(p,F_n)^*) \]
	where $r$ runs over the positive divisors of $f$.
	We compute separately the indices $(R:U_f)$ and $(U_f:U)$.
	
	Recall that if $p|f$ then $T_p=\Delta$, where $\Delta=\gal(F_n|\Q_n)$. We have
	\[ s(\Delta)(1-(p,F_n)^*)=s(\Delta)(1-\frac{s(\Delta)}{2}\lambda_p^{-1})=s(\Delta)(1-\lambda_p^{-1})\in R. \]
	Moreover, choosing $r=f$ we get that $s(\Delta)\in U_f$. Therefore, we deduce that
	\[ U_f=s(\Delta)R+\prod_{p|f}(1-(p,F_n)^*)R. \]
	Denote by $e_{\Delta}$ the element $s(\Delta)/2$. For any $p$ that  divides $f$ we have
	\[ (1-e_{\Delta})(1-(p,F_n)^*)=(1-e_{\Delta}).\]
	Moreover, $(1-e_{\Delta})s(\Delta)=0$. Let $x\in U_f$, there exist $u,v\in R$ such that the element $x$ is equal to $s(\Delta)u+\prod_{p|f}(1-(p,F_n)^*)v$. Then $(1-e_{\Delta})x=(1-e_{\Delta})v$ and hence
	\[ (1-e_{\Delta})U_f=(1-e_{\Delta})R. \]
	Note that for any $R$-module $M\subset\Q(\Gamma_n)$ the kernel of the map $M\longrightarrow (1-e_{\Delta})M$ induced by the multiplication by $1-e_{\Delta}$ is exactly $M^{\Delta}$.
	
	We show that $U_f^{\Delta}=R^{\Delta}$. Indeed,	it is clear that $R^{\Delta}=s(\Delta)R$. Then, take $x\in U_f^{\Delta}$. There exist $u,v\in R$ such that $x=s(\Delta)u+\prod_{p|f}(1-(p,F_n)^*)v$ and $(1-e_{\Delta})v=0$. Therefore, the element $v$ is divisible by $s(\Delta)$. It follows that
	\[ U_f^{\Delta}=s(\Delta)R+\prod_{p|f}(1-(p,F_n)^*)s(\Delta)R=s(\Delta)R+\prod_{p|f}(1-\lambda_p^{-1})s(\Delta)R=s(\Delta)R. \] 
	
	Finally, from the Equation \ref{index2} we deduce
	\[ (R:U_f)=(R^{\Delta}:U_f^{\Delta})((1-e_{\Delta})R:(1-e_{\Delta})U_f)=1. \]
	
	We compute now the index $(U_f:U)$. Let $U_f(T_2)$ be the $R$-submodule of $U$ generated by
	\[ s(T_{2r})\prod_{p|\frac{2f}{2r}}(1-(p,F_n)^*) \]
	where $r$ runs over the positive divisors of $f$. Then, we get
	\[ U=U_f(T_2)+(1-(2,F_n)^*)U_f. \]
	By \cite[Lemma $5.1$]{sinnott} the index $(U_f:U)$ is equal to the cardinality of the quotient
	\begin{equation}\label{Sinnott-index-quotient} \tag*{($A3$)}
		\frac{U_f^{T_2}}{U_f(T_2)+(1-\lambda_2^{-1})U_f^{T_2}}.
	\end{equation} 
	Notice that the inertia group $T_2$ is $\Gamma_n$ if $f\equiv1\bmod4$ or $T_2=\gal(F_n|F)$ if $f\equiv3\bmod4$. Then, it is easy to  show that
	\[ U_f(T_2)=s(\Gamma_n)R+s(T_2)\prod_{p|f}(1-(p,F_n)^*)R. \]
	
	Recall that $U_f^{\Delta}=R^{\Delta}$ and $(1-e_{\Delta})U_f=(1-e_{\Delta})R$ as shown previously. Consider the following commutative diagram whose rows are exact
	\[ \xymatrix{0\ar[r]&R^{\Delta}\ar[d]^-{\lambda'}\ar[r]&R\ar[d]^-{\lambda}\ar[r]^-{1-e_{\Delta}}&(1-e_{\Delta})R\ar[r]\ar@{=}[d]&0\\
		0\ar[r]&R^{\Delta}\ar[r]&U_f\ar[r]^-{1-e_{\Delta}}&(1-e_{\Delta})R\ar[r]&0} \]
	The application $\lambda'$ is the multiplication by $\prod_{p|f}(1-\lambda_p^{-1})$ while the application $\lambda$ is the multiplication by $\prod_{p|f}(1-(p,F_n)^*)$. The commutativity of the diagram follows from the relations
	\begin{align*} s(\Delta)(1-(p,F_n)^*)=s(\Delta)(1-\lambda_p^{-1}) \\ (1-e_{\Delta})(1-(p,F_n)^*)=1-e_{\Delta} 
	\end{align*}
	for any prime divisor $p$ of $f$. We take $T_2$-invariants in order to get the commutative diagram
	\[ \xymatrix{0\ar[r]&R^{\Gamma_n}\ar[d]^-{\lambda'}\ar[r]&R^{T_2}\ar[d]^-{\lambda}\ar[r]^-{1-e_{\Delta}}&((1-e_{\Delta})R)^{T_2}\ar[r]\ar@{=}[d]&0\\
		0\ar[r]&R^{\Gamma_n}\ar[r]&U_f^{T_2}\ar[r]^-{1-e_{\Delta}}&((1-e_{\Delta})R)^{T_2}\ar[r]&0}. \]
	By \cite[Lemma $5.2$]{sinnott} we have that $\mathrm{H}^1(T_2,R^{\Delta})=0$. It follows the exactness of the diagram's rows. We deduce from the diagram that
	\[ U_f^{T_2}=R^{\Gamma_n}+\lambda\left( R^{T_2}\right)=s(\Gamma_n)R+s(T_2)\prod_{p|f}(1-(p,F_n)^*)R.  \]
	In other words, we get that $U_f^{T_2}=U_f(T_2)$. Then, the quotient appearing in \ref{Sinnott-index-quotient} is trivial and the index $(U_f:U)$ is equals to $1$.
	
	By Equation \ref{index1}, we conclude that
	\[ c_{F_n}=(R:U)=(R:U_f)(U_f:U)=1. \]
	
\end{proof}

\section{More cyclotomic units}\label{appB}

In this appendix we collect some information about the units defined in Definition \ref{moreunits} and we prove Theorem \ref{CFn-structure}. In particular, we build the group $C_{F_n}$ by taking the possible square roots of Sinnott cyclotomic units.

We summarize some useful properties of the units appearing in Definition \ref{moreunits}. Before, we define the unit $\delta_f$ for any squarefree integer $f$.
\begin{definition}\label{full-deltaf}
	Fix an algebraic closure $\overline{\Q}$ of the rationals. For any positive squarefree integer $f$, let $F=\Q(\sqrt{f})$ and denote by $\sigma$ the generator of $\Delta=\gal(F|\Q)$ . Let $\delta_f\in\overline{\Q}$ be such that
	\[ \delta_f^2=\begin{cases}
		\underset{\Q(\ze{4f})|F}{\norm}(1-\ze{4f}), & \t{if $f\equiv2,3\bmod4$,}\\
		\left(\underset{\Q(\ze{f})|F}{\norm}(1-\ze{f})\right)^{1-\sigma}, & \t{if $f\equiv1\bmod4$ and $f$ prime,}\\
		\underset{\Q(\ze{f})|F}{\norm}(1-\ze{f}), & \t{otherwise.}
	\end{cases} \]
	Note that $\delta_f$ is defined up to a sign and $\delta_f^2$ lies in $O_{F_0}^{\times}$.
\end{definition}

\begin{lemma}\label{square-properties}
	The units in Definition \emph{\ref{moreunits}} have the following properties:
	\begin{enumerate}[label=(\alph*)]
		\item\label{prop-beta} Let $\gamma$ be a generator of $G_n$.  The square of $\beta_n$ is a cyclotomic unit of $\Q_n$, namely
		\[ \beta_n^2=\left(\underset{\Q(\ze{2^{n+2}})|\Q_n}{\norm}(1-\ze{2^{n+2}})\right)^{\gamma-1}. \]		
		\item\label{prop-eta} The square of $\eta_n$ is a cyclotomic unit of $F_n$, namely
		\[ \eta_n^2=\underset{\Q(\ze{2^{n+2}f})|F_n}{\norm}(1-\ze{2^{n+2}f}). \]
		\item\label{prop-deltaf} Assume $f\equiv1\bmod4$ and $f$ prime. Let $s\in(\Z/f\Z)^{\times}$ be such that the automorphism of $\Q(\ze{f})$ defined by $\ze{f}\mapsto\ze{f}^s$ generates $\Delta=\gal(F|\Q)$ when restricted to the field $F$. Then up to a sign we get
		\[ \delta_f=\underset{\Q(\ze{f})^+|F}{\norm}\left(\ze{f}^{\frac{s-1}{2}}\frac{1-\ze{f}}{1-\ze{f}^s}\right). \]
	\end{enumerate}
	
\end{lemma}

\begin{proof}
	All the three parts are simple computations. They follow from the relations:
	\begin{align*}
		\left( \ze{2^{n+2}}^{\frac{1-c}{2}}\frac{1-\ze{2^{n+2}}^c}{1-\ze{2^{n+2}}}\right)^2  =& \frac{(1-\ze{2^{n+2}}^c)(1-\ze{2^{n+2}}^{-c})}{(1-\ze{2^{n+2}})(1-\ze{2^{n+2}}^{-1})},\\
		\left( \ze{4}(\ze{2^{n+3}f}-\ze{2^{n+3}f}^{-1})\right)^2  =& (1-\ze{2^{n+2}f})(1-\ze{2^{n+2}f}^{-1}),\\
		\left( \ze{f}^{\frac{s-1}{2}}\frac{1-\ze{f}}{1-\ze{f}}^s\right)^2  =& \frac{(1-\ze{f})(1-\ze{f}^{-1})}{(1-\ze{f}^s)(1-\ze{}^{-s})}.
	\end{align*}
	
	For part \ref{prop-deltaf}, we recall that if $f\equiv1\bmod4$ and $f$ is prime then we have
	\[ \delta_f^2=\left(\underset{\Q(\ze{f})|F}{\norm}(1-\ze{f})\right)^{1-\sigma}=\underset{\Q(\ze{f})|F}{\norm}\left(\frac{1-\ze{f}}{1-\ze{f}^s}\right). \]
	
\end{proof}

\begin{lemma}\label{cyclotomic-units}
	Let $f\ge3$ be an odd squarefree integer. Then, for any $n\ge0$ the group of cyclotomic units for the field $F_n$ is
	\[ \cyc_{F_n}=\left\langle -1,\eta_n^2,\beta_n^2,\delta_f^2\right\rangle_{\Z[G_n]} .\]
\end{lemma}

\begin{proof}
	We prove the statement for $n\ge1$. The case $n=0$ is similar and easier. Define the following units:
	\begin{align*}
		a_t=&\underset{\Q(\ze{2^{n+2}f})|F_n}{\norm}(1-\ze{2^{n+2}f}),\quad 0\le t\le n,\\
		b_t=&\underset{\Q(\ze{2^{n+2}})|\Q_n}{\norm}(1-\ze{2^{n+2}}),\quad 1\le t\le n,\\
		d_f=&\underset{\Q(\ze{f'})|F}{\norm}(1-\ze{f'}).
	\end{align*}
	where $f'$ is the conductor of $F$, namely $f'=f$ if $f\equiv1\bmod4$ or $f'=4f$ if $f\equiv3\bmod4$. Notice that $\norm_{F_n|F_t}(a_n)=a_t$ and $\norm_{\Q_n|\Q_t}(b_n)=b_t$. Combining this information with Remark \ref{cyclo-DK} we get that
	\[ D_{F_n}=\Q^\times\cdot\left\langle a_n,b_n,d_f\right\rangle_{\Z[\Gamma_n]}.  \]
	Moreover, we can consider only the $G_n$-conjugates. Indeed, let $\sigma$ be the generator of $\Delta$, then $a_n^{1+\sigma}$ lies in $\left\langle -1,b_n\right\rangle_{\Z[G_n]} $, the element $b_n$ is $\Delta$-invariant and $d_f^{1+\sigma}\in\Q^\times$. Let $\gamma$ be a generator of $G_n$, then 
	\[ D_{F_n}=\Q^\times\cdot\left\langle  a_n^{\gamma^i}, b_n^{\gamma^j}, d_f\;:\; 0\le i < 2^n,\; 0\le j <2^n\right\rangle_{\Z}.  \]
	By Lemma \ref{square-properties} we have $a_n=\eta_n^2$ and $b_n^{\gamma-1}=\beta_n^2$. Therefore
	\[ D_{F_n}=\Q^\times\cdot\left\langle  (\eta_n^2)^{\gamma^i}, b_n, (\beta_n^2)^{\gamma^j}, d_f\;:\; 0\le i < 2^n,\; 0\le j <2^n-1\right\rangle_{\Z}.  \]
	
	According to \cite[Lemma $4.1$]{sinnott}, an element $u\in D_{F_n}$ lies in $\cyc_{F_n}$ if and only if $\norm_{F_n|\Q}(u)=\pm1$. We compute the norm of the generators:
	\begin{align*}
		&\underset{F_n|\Q}{\norm}((\eta_n^2)^{\gamma^i})=1,\\
		&\underset{F_n|\Q}{\norm}((\beta_n^2)^{\gamma^j})=1,\\
		&\underset{F_n|\Q}{\norm}(b_n)=4,\\
		&\underset{F_n|\Q}{\norm}(d_f)=\begin{cases}f^{2^n}, & \t{if $f\equiv1\bmod4$ and $f$ prime,}\\
			1, & \t{otherwise.}\end{cases}
	\end{align*}
	
	Assume now that $f\equiv1\bmod4$ and $f$ prime. The elements with norm $\pm1$ in $D_{F_n}$ are
	\[ \cyc_{F_n}=\left\langle -1, (\eta_n^2)^{\gamma^i}, (\beta_n^2)^{\gamma^j}, b_n^{2^n}/2, d_f^2/f\right\rangle_{\Z}.  \]
	Note that
	\[ \frac{d_f^2}{f}=\frac{d_f^{1+\sigma}d_f^{1-\sigma}}{f}=d_f^{1-\sigma}=\delta_f^{2}. \]
	Consider now the element $b_n^{2^n}/2$. Let $N_n=\sum_{i=0}^{2^n-1}\gamma^i\in\Z[G_n]$, then $2^n-N_n$ lies in the augmentation ideal. In particular, we have $2^n-N_n=(\gamma-1)x$ for some $x\in\Z[G_n]$. So
	\[ \frac{b_n^{2^n}}{2}=\frac{b_n^{N_n}b_n^{2^n-N_n}}{2}=b_n^{(\gamma-1)x}=(\beta_n^2)^x. \]
	We conclude that
	\[ \cyc_{F_n}=\left\langle -1,\eta_n^2,\beta_n^2,\delta_f^2\right\rangle_{\Z[G_n]}.  \]
	
	In the remaining cases $d_f$ has norm $1$. Therefore, we get that
	\[ \cyc_{F_n}=\left\langle -1, (\eta_n^2)^{\gamma^i}, (\beta_n^2)^{\gamma^j}, b_n^{2^n}/2, d_f\right\rangle_{\Z}.  \]
	By Definition \ref{full-deltaf} we have that $d_f=\delta_f^2$. Moreover we can show as above that $b_n^{2^n}/2$ lies in the $\Z[G_n]$-submodule generated by $\beta_n^2$. We finally get
	\[  \cyc_{F_n}=\left\langle -1,\eta_n^2,\beta_n^2,\delta_f^2\right\rangle_{\Z[G_n]}.  \]  
\end{proof}

\begin{lemma}\label{eta-is-unit}
	Let $f\ge3$ be an odd squarefree integer. Assume that the class number of $\Q(\sqrt{f})$ is even if $f\not\equiv1\bmod8$. Then $\eta_n\in O_{F_n}^{\times}$.
\end{lemma}

\begin{proof}
	We compute the norm of $\eta_n^2$. First, we assume $f\equiv3\bmod4$. We get
	\begin{equation}\label{norm34}\tag*{($B1$)}
		\underset{F_n|F}{\norm}(\eta_n^2)=\underset{\Q(\ze{2^{n+2}f})|F}{\norm}(1-\ze{2^{n+2}f})=\underset{\Q(\ze{4f})|F}{\norm}(1-\ze{4f})=\eta_0^2=\delta_f^2.
	\end{equation}
	
	By Sinnott's index formula \ref{index-sinnott}, we have that the index $[O_{F}^{\times}:\cyc_F]=\#\cl_F$ and it is even. Since $\cyc_F=\left\langle -1,\delta_f^2\right\rangle $ and $\delta_f^2$ is positive under any real embedding of $F$ we deduce that $\delta_f^2$ is a square in $O_{F}^{\times}$.
	
	Assume now that $f\equiv1\bmod4$. Let $\phi_2$ be the automorphism of $\Q(\ze{2^{n+2}f})$ such that $\phi_{2|\Q(\ze{2^{n+2}})}=\t{id}$ and $\phi_2(\zeta_f)=\zeta_f^2$. Then
	\begin{align*}
		\underset{F_n|F}{\norm}(\eta_n^2)=&\underset{\Q(\ze{2^{n+2}f})|F}{\norm}(1-\ze{2^{n+2}f})=\\&=\underset{\Q(\ze{f})|F}{\norm}\left(({\tiny } 1-\ze{f})^{1-\phi_2^{-1}}\right)=\left(\underset{\Q(\ze{f})|F}{\norm}\left( 1-\ze{f}\right)\right)^{1-\phi_{2|F}}.  
	\end{align*}
	Let $\chi_f$ be the quadratic Dirichlet character associated to $F$. Then $\phi_{2|F}=\t{id}$ if and only if $2\in\ker\chi_f$. The last condition holds if and only if $2$ splits in ${F}$, namely when $f\equiv1\bmod8$. Therefore we get
	\begin{equation}\label{norm14}\tag*{($B2$)}
		\underset{F_n|F_0}{\norm}(\eta_n^2)=\begin{cases}
			1, & \t{if $f\equiv1\bmod8$},\\
			\delta_f^2, & \t{if $f\equiv5\bmod8$ and $f$ prime,}\\
			\delta_f^{2(1-\sigma)}=\delta_f^4, & \t{otherwise.}
		\end{cases}
	\end{equation} 
	Notice that the right hand side is always a square in $O_{F}^\times$. See Lemma \ref{square-properties} part \ref{prop-deltaf} for the case $f\equiv5\bmod8$ and $f$ prime.
	
	Let $\gamma$ be a generator of the group $G_{n+1}$.	Let $N_n^{(n+1)}=\sum_{i=0}^{2^n-1}\gamma^i\in\Z[G_{n+1}]$. We have shown that $\norm_{F_n|F}(\eta_n^2)=\eta_n^{2N_n^{(n+1)}}$ is a square in $O_{F}^\times$. It follows that $\eta_n^{N_n^{(n+1)}}\in O_{F}^\times$ and we get that
	\[ \eta_n^{\gamma^{2^n}-1}=\left( \eta_n^{N_n^{(n+1)}}\right)^{\gamma-1}=1.  \]
	Since $\gamma^{2^n}$ generates $\gal(F_{n+1}|F_n)$, we have $\eta_n\in O_{F_n}^{\times}$.
	
\end{proof}

\begin{remark}\label{norm-relation}
	We have the following relation between the units $\eta_n$:
	\[ \underset{F_n|F_m}{\norm}(\eta_n^2)=\underset{\Q(\ze{2^{n+2}f})|F_m}{\norm}(1-\ze{2^{n+2}f})=\underset{\Q(\ze{2^{m+2}f})|F_m}{\norm}(1-\ze{2^{m+2}f})=\eta_m^2. \]
	In particular, taking the square roots in $O_{F_n}^{\times}$ we get that $\norm_{F_n|F_m}\eta_n=\eta_m$ up to a sign. Moreover, we have seen in equations \ref{norm34} and \ref{norm14} that up to a sign
	\[ \eta_0=\begin{cases}
		1, & \t{if $f\equiv1\bmod8$, }\\
		\delta_f^2, & \t{if $f\equiv5\bmod8$ and $f$ is not prime,}\\
		\delta_f, & \t{otherwise}.
	\end{cases} \]
	
	Similar relations can be deduced for $\beta_n$. Indeed, letting $\gamma$ be a generator of $G_n$ we have
	\[ \underset{\Q_n|\Q_m}{\norm}(\beta_n^2)=\left(\underset{\Q(\ze{2^{n+2}})|\Q_m}{\norm}(1-\ze{2^{n+2}})\right)^{\gamma-1}=\left(\underset{\Q(\ze{2^{m+2}})|\Q_m}{\norm}(1-\ze{2^{m+2}})\right)^{\gamma-1}=\beta_m^2.  \]
	We deduce that $\norm_{\Q_n|\Q_m}\beta_n=\beta_m$ up to a sign.
	
\end{remark}

We prove now Theorem \ref{CFn-structure}.
\begin{proof}[Proof of Theorem \emph{\ref{CFn-structure}}.]
	First of all, combining Lemmas \ref{square-properties}, \ref{cyclotomic-units} and \ref{eta-is-unit} we get that
	\[ \cyc_{F_n}\subset C_{F_n}\subset O_{F_n}^{\times}. \]
	Since $\cyc_{F_n}$ has finite index in $O_{F_n}^{\times}$ as follows from Sinnott's index formula, the module $C_{F_n}$ and $O_{F_n}^{\times}$ have the same $\Z$-rank. Moreover, the only torsion elements are $\pm1$. So, the morphisms $\pi$ of Theorem $1.4$ are surjective morphisms between free $\Z$-modules of the same rank, therefore they are isomorphisms.
	
	In order to prove the statement about the index, we notice that the index $[C_{F_n}:\cyc_{F_n}]$ is equal to $[\pi^{-1}(C_{F_n}/\pm1):\pi^{-1}(\cyc_{F_n}/\pm1)]$. Furthermore, by Lemma \ref{cyclotomic-units} we have
	\begin{align*}
		\pi^{-1}(\cyc_{F_n}/\pm1)=\left<2\right>_{\Z[G_n]}\oplus\left<2\right>_{\Z[G_n]/(N_n)}, & \quad\t{if $f\equiv3\bmod4$,}\\
		\pi^{-1}(\cyc_{F_n}/\pm1)=\left<2,2^jN_n\right>_{\Z[G_n]}\oplus\left<2\right>_{\Z[G_n]/(N_n)}, & \quad\t{if $f\equiv5\bmod8$,}\\
		\pi^{-1}(\cyc_{F_n}/\pm1)=\left<2\right>_{\Z[G_n]}\oplus\left<2\right>_{\Z[G_n]/(N_n)}\oplus\left<2^j\right>_{\Z}, & \quad\t{if $f\equiv1\bmod8$,}
	\end{align*} 
	where $j=0$ if $f$ is prime $j=1$ otherwise. So, we compute the index
	\[ [\pi^{-1}(C_{F_n}/\pm1):\pi^{-1}(\cyc_{F_n}/\pm1)]=\begin{cases}
		2^{2^{n+1}-1}, & \t{if $f\equiv1\bmod4$ and $f$ is not prime,}\\
		2^{2^{n+1}-2}, & \t{otherwise.}
	\end{cases} \]
	Then, combining with Sinnott's index formula $\ref{index-sinnott}$ and the fact that $c_{F_n}=1$ we get $[O_{F_n}^{\times}:C_{F_n}]=\#\cl_{F_n}$ if $f\equiv1\bmod4$ or $[O_{F_n}^{\times}:C_{F_n}]=\#\cl_{F_n}/2$ if $f\equiv3\bmod4$. 
	
	We need to prove that $C_{F_n}^{\gal(F_n|F_m)}=C_{F_m}$. Assume $f\not\equiv1\bmod8$. Under the isomorphism $\pi$ we get that
	\[ \underset{F_n|F_m}{\norm}\left(\Z[G_n]\oplus\frac{\Z[G_n]}{(N_n)}\right)\overset{\pi}{\cong}\left(\frac{C_{F_n}}{\pm1}\right)^{\gal(F_n|F_m)}. \]
	We deduce that $(C_{F_n}/\pm1)^{\gal(F_n|F_m)}$ is generated as a $\Z[G_n]$-modules by the elements $\norm_{F_n|F_m}\eta_n=\pm\eta_m$ and $\norm_{F_n|F_m}\beta_n=\pm\beta_m$. These norms are computed in Remark \ref{norm-relation}. We have a natural injective map
	\[ \frac{C_{F_n}^{\gal(F_n|F_m)}}{\pm1}\longrightarrow \left(\frac{C_{F_n}}{\pm1}\right)^{\gal(F_n|F_m)} \]
	and since both $\eta_m$ and $\beta_m$ are invariant under the action of $\gal(F_n|F_m)$ we deduce that
	\[ C_{F_n}^{\gal(F_n|F_m)}=\left<-1,\eta_m,\beta_m\right>=C_{F_m}. \]
	
	Assume now $f\equiv1\bmod8$, we have
	\[ \underset{F_n|F_m}{\norm}\left(\frac{\Z[G_n]}{(N_n)}\oplus\frac{\Z[G_n]}{(N_n)}\right)\oplus\Z\overset{\pi}{\cong}\left(\frac{C_{F_n}}{\pm1}\right)^{\gal(F_n|F_m)}. \]
	Similarly as above, it is easy to deduce that 
	\[ C_{F_n}^{\gal(F_n|F_m)}=\left<-1,\eta_m,\beta_m,\delta_f'\right>=C_{F_m}. \]
\end{proof}

We conclude proving some lemmas about these units.

\begin{lemma}\label{naka2}
	Let $A,B$ be free $\Z$-modules such that $A\subset B$. Then, the morphisms induced by the inclusion
	\[ \iota_k:A/2^kA\longrightarrow B/2^kB \]
	are injective for each $k\ge1$ if and only if $\iota_1$ is injective.
\end{lemma}

\begin{lemma}\label{Qn-injective}
	Let $C_{\Q_n}$ be the $\Z[G_n]$-submodule of $O_{\Q_n}^{\times}$ generated by $-1$ and $\beta_n$. Then, the natural morphisms
	\[ \iota:\frac{C_{\Q_n}}{C_{\Q_n}^{2^k}}\longrightarrow \frac{O_{F_n}^\times}{O_{F_n}^{\times2^k}}\]
	are injective for any $n,k\ge1$.
\end{lemma}
\begin{proof}
	The kernel of $\iota$ is isomorphic to the kernel of
	\[ \frac{C_{\Q_n}}{\pm C_{\Q_n}^{2^k}}\longrightarrow \frac{O_{F_n}^\times}{\pm O_{F_n}^{\times2^k}}. \]
	Since $\Q_n$ and $F_n$ are both real fields, we can apply Lemma \ref{naka2} to the free modules $C_{\Q_n}/\pm1$ and $O_{F_n}^{\times}/\pm1$ so that it suffices to prove the injectivity of
	\[ \iota':\frac{C_{\Q_n}}{\pm C_{\Q_n}^{2}}\longrightarrow \frac{O_{F_n}^\times}{\pm O_{F_n}^{\times2}}. \]

	By \cite[Theorem $8.2$]{wash} the index $[O_{\Q_n}^{\times}:C_{\Q_n}]$ is equal to the class number of $\Q_n$ and it is odd. Therefore, the natural map 
	\[ C_{\Q_n}/\pm C_{\Q_n}^{2}\cong O_{\Q_n}^\times/\pm O_{\Q_n}^{\times 2} \]
	is an isomorphism.
	
	Since $-1$  is not a square in $O_{F_n}^{\times}$, the kernel of $\iota'$ is isomorphic to the kernel of the map 
	\[ \frac{O_{\Q_n}^{\times}}{O_{\Q_n}^{\times2}}\longrightarrow\frac{O_{F_n}^{\times}}{O_{F_n}^{\times2}}. \]
	
	Assume that there exists $u\in O_{\Q_n}^{\times}$ that is a square in ${F_n}^{\times}$ but not in $\Q_n^{\times}$. Then, we have that 
	\[ F_n=\Q_n(\sqrt{u}). \]
	According to \cite[Theorem $6.3$]{gras} the extension $\Q_n\subset\Q_n(\sqrt{u})$ can only be ramified at the prime above $2$. This leads to a contradiction since $F_n=\Q_n(\sqrt{f})$ with $f$ an odd squarefree integer.  
\end{proof} 

\begin{lemma}\label{Fn-Qn-injectivity}
	Assume that the prime $2$ does not ramify in $F$. Let $C_{\Q_n}$ be the $\Z[G_n]$-submodule of $C_{F_n}$ generated by $-1$ and $\beta_n$. Then, the natural morphisms
	\[ \frac{C_{\Q_n}}{\pm C_{\Q_n}^{2^k}}\oplus\frac{O_{F}^\times}{\pm O_{F}^{\times 2^k}}\longrightarrow\frac{O_{F_n}^\times}{\pm O_{F_n}^{\times 2^k}} \]
	are injective for any $n,k\ge1$.
\end{lemma}
\begin{proof}
	As a consequence of Lemma \ref{naka2}, it suffices to show that the map
	\[ \iota:\frac{C_{\Q_n}}{\pm C_{\Q_n}^{2}}\oplus\frac{O_{F}^\times}{\pm O_{F}^{\times 2}}\longrightarrow\frac{O_{F_n}^\times}{\pm O_{F_n}^{\times 2}} \]
	is injective. 
	
	Consider the following commutative diagram of $\Z[G_n]$-modules:
	\[ \xymatrix{0\ar[r]&0\ar[r]\ar[d]&\frac{C_{\Q_n}}{\pm C_{\Q_n}^2}\ar[r]\ar[d]&\frac{C_{\Q_n}}{\pm C_{\Q_n}^2}\ar[r]\ar[d]&0\\
		0\ar[r]&\ker\iota\ar[r]&\frac{C_{\Q_n}}{\pm C_{\Q_n}^{2}}\oplus\frac{O_{F}^\times}{\pm O_{F}^{\times 2}}\ar[r]&\frac{O_{F_n}^\times}{\pm O_{F_n}^{\times 2}}   }.\]
	The rows are exact. Moreover, by Lemma \ref{Qn-injective} the rightmost vertical map is injective. Using the snake lemma we can identify $\ker\iota$ with a subgroup of $O_{F}^{\times}/\pm O_{F}^{\times2}$. Therefore, the cardinality of $\ker\iota$ is at most $2$. 
	
	We know that
	\[ \frac{C_{\Q_n}}{\pm C_{\Q_n}^2}\oplus\frac{O_{F}^{\times}}{\pm O_{F}^{\times2}}\cong\frac{\Z[G_n]}{(2,N_n)}\oplus\frac{\Z}{2\Z}. \]
	So the only $\Z[G_n]$-submodules with $2$ elements are the following:
	\[ \left<(\beta_1,1)\right>,\quad \left<(1,\varepsilon)\right>,\quad \left<(\beta_1,\varepsilon)\right>, \]
	where $\varepsilon$ is a fundamental unit of $F$.
	Note that $\beta_1\not\in \pm O_{F_n}^{\times^2}$ by Lemma \ref{Qn-injective}.
	
	Assume that $\ker\iota$ is generated by $(1,\varepsilon)$. Then, the unit $\varepsilon$ lies in $\pm O_{F_n}^{\times2}$. Since the fundamental unit $\varepsilon$ is defined up to a sign, we assume that there exists $z\in O_{F_n}^{\times}$ such that $\varepsilon=z^2$. Otherwise, if $\varepsilon=-z^2$ we consider $-\varepsilon$ as a fundamental unit. Note that the extension $F(z)$ has degree $2$ over $F$.  Therefore, $z\in F_1=F(\sqrt{2})$. So $z=x+\sqrt{2}y$ for some $x,y\in F$. The relation $z^2=\varepsilon$ gives us two cases: or $x^2=\varepsilon$ or $2y^2=\varepsilon$. The first case is impossible as $\varepsilon$ is a fundamental unit. On the other hand, if there exists $y\in F_0$ such that $2y^2=\varepsilon$, then $(y^{-1})^2=2\varepsilon^{-1}$. In terms of ideals, we get that 
	\[ 2O_{F}=\left(y^{-1}O_{F}\right)^2\]
	This leads to a contradiction since the prime $2$ does not ramify in $F$.
	
	Assume that $\ker\iota$ is generated by $(\beta_1,\varepsilon)$. Then, the unit $\beta_1\varepsilon$ lies in $\pm O_{F_n}^{\times2}$. So there exists $z\in O_{F_n}^{\times}$ such that
	\[ \beta_1\varepsilon=z^2, \quad\t{or}\quad \beta_1\varepsilon=-z^2. \]
	Consider a generator $\gamma$ of $G_n=\gal(F_n|F)$. Then
	\[ ( z^{1+\gamma})^2=(\pm z^2)^{1+\gamma}=\underset{\Q_1|\Q}{\norm}(\beta_1)\varepsilon^2=-\varepsilon^2. \] 
	This leads to a contradiction with the fact that $F_n$ is a real field.
	
	We conclude that $\ker\iota=1$.
\end{proof}

\vspace{1.5cm}
%\nocite{}


\begin{thebibliography}{99}
		
\bibitem[Che33]{chevalley} Chevalley C.: \emph{Sur la théorie du corps de classes dans les corps finis et les corps locaux}, Jour. of the Fac. of Sc., Tokyo, Vol. II, (1933).

\bibitem[Cor97]{corn} Cornacchia P.: \emph{Anderson's module for cyclotomic fields of prime conductor}, J. Number Theory \textbf{67}(2) (1997), 252-276.

\bibitem[FW79]{ferrerowash} Ferrero B. and Washington L.C.: \emph{The Iwasawa invariant $\mu_p$ vanishes for abelian number fields}, Ann. of Math. Second Series, Vol. \textbf{109}, No. 2 (1979), 377–395.

\bibitem[Fuk94]{fukuda-stab} Fukuda T.: \emph{Remarks on $\Z_p$-extensions of number fields}, Proc. Japan Acad. Ser. A \textbf{70}(8) (1994), 264–266.

\bibitem[Fuk96]{fukuda} Fukuda T.: \emph{Cyclotomic units and Greenberg's conjecture for real quadratic fields}, Math. Comp. Volume \textbf{65}, Number 215 (1996), 1339-1348.

\bibitem[FK05]{fukudakomatsuprev} Fukuda T. and Komatsu K.: \emph{On the Iwasawa $\lambda$-invariant of the cyclotomic $\Z_2$-extension of a real quadratic field}, Tokyo J. Math. \textbf{28}(1) (2005), 259-264.

\bibitem[FK09]{fukudakomatsu} Fukuda T. and Komatsu K.: \emph{On the Iwasawa $\lambda$-invariant of the cyclotomic $\Z_2$-extension of $\Q(\sqrt{p})$}, Math. Comp. Vol. \textbf{78}, No. 267 (2009), 1797-1808.

\bibitem[Gra05]{gras} Gras G.: \emph{Class Field Theory: from theory to practice} , corr. 2nd ed., SMM (2005).

\bibitem[Gre76]{green} Greenberg, R.: \emph{On the Iwasawa invariants of totally real number fields}, Amer. J. Math. \textbf{98}(1) (1976), 263-284.

\bibitem[KS95]{kraftschoof} Kraft J. and Schoof R.: \emph{Computing Iwasawa modules of real quadratic number fields}, Compos. Math. \textbf{97} Issues 1-2 (1995), 135-155 + Erratum. Compos. Math., \textbf{103} no. 2 (1996), p. 241-241 

\bibitem[Iwa59]{iwasawa} Iwasawa K.: \emph{On $\Gamma$-extensions of algebraic number fields}, Bull. Amer. Math. Soc. \textbf{65}(4) (1959), 183–226.

\bibitem[Mer21]{mercuri} Mercuri P.: unpublished computations, 2021.

\bibitem[Miz04]{mizusawa} Mizusawa Y.: \emph{On the Iwasawa invariants of $\Z_2$-extensions of certain real quadratic fields}, Tokyo J. Math. \textbf{27}(1) (2004), 255-261.

\bibitem[Nis06]{nishino} Nishino Y.: \emph{On the Iwasawa invariants of the cyclotomic $\Z_2$-extension of certain real quadratic fields}, Tokyo J. Math. \textbf{29}(1) (2006), 239-245.

\bibitem[OT97]{ozakitaya} Ozaki M. and Taya H.:, \emph{On the Iwasawa $\lambda_2$-invariants of certain family of real quadratic fields}, Manuscripta
Math. \textbf{94} (1997), 437–444.

\bibitem[Pao02]{paoluzi} Paoluzi M.: \emph{La congettura di Greenberg per campi quadratic reali}, Ph.D. dissertation at Università di Roma ``La Sapienza'', 2002. 

\bibitem[Sch02]{schoof2} Schoof R.: \emph{Class numbers of real cyclotomic fields of prime conductor}, Math. Comp. Volume \textbf{72} Number 242 (2002), 913-937.

\bibitem[Sin80]{sinnott} Sinnott W.: \emph{On the Stickelberger ideal and the circular units of an abelian field}, Invent. Math. \textbf{62} (1980), 181-234.

\bibitem[Tay96]{taya} Taya H.: \emph{Computation of $\Z_3$-invariants of real quadratic fields}, Math. Comp. Volume \textbf{65} Issue 214 (1996), 779-784.

\bibitem[Was97]{wash} Washington L.C.: \emph{Introduction to cyclotomic fields}, Second Edition, Graduate Texts in Math. \textbf{83}, Springer-Verlag, Berlin Heidelberg New York, (1997). 

		
\end{thebibliography}
\end{document}